\pdfoutput=1

\documentclass[preprint,10pt]{elsarticle}

\usepackage{fullpage}
\usepackage{hyperref}
\usepackage{amsmath,amssymb,amsfonts,mathrsfs}
\usepackage[titletoc,toc,title]{appendix}

\usepackage{array} 
\usepackage[utf8]{inputenc}
\usepackage{listings}
\usepackage{mathtools}
\usepackage{pdfpages}
\usepackage[textsize=footnotesize,color=green]{todonotes}
\usepackage{bm}
\usepackage[normalem]{ulem}
\usepackage{hhline}

\usepackage{algorithm}
\usepackage[noend]{algpseudocode}
\usepackage{algorithmicx}
\algblock{ParFor}{EndParFor}
\algnewcommand\algorithmicparfor{\textbf{parfor}}
\algnewcommand\algorithmicpardo{\textbf{do}}
\algnewcommand\algorithmicendparfor{\textbf{end\ parfor}}
\algrenewtext{ParFor}[1]{\algorithmicparfor\ #1\ \algorithmicpardo}
\algrenewtext{EndParFor}{\algorithmicendparfor}

\usepackage{graphicx}
\usepackage{subfig}
\usepackage{color}

\usepackage{pgfplots}
\usepackage{pgfplotstable}
\definecolor{markercolor}{RGB}{124.9, 255, 160.65}
\pgfplotsset{width=10cm,compat=1.9}
\pgfplotsset{
tick label style={font=\small},
label style={font=\small},
legend style={font=\small}
}

\usetikzlibrary{calc}
\usepackage{stmaryrd}

\newcommand{\pd}[2]{\frac{\partial#1}{\partial#2}}

\newcommand{\mi}{\mathrm{i}} 

\newcommand{\nor}[1]{\left\| #1 \right\|}

\newcommand{\LRp}[1]{\left( #1 \right)}

\newcommand{\LRa}[1]{\left\langle #1 \right\rangle}
\newcommand{\LRb}[1]{\left| #1 \right|}
\newcommand{\LRc}[1]{\left\{ #1 \right\}}

\newcommand{\Grad} {\ensuremath{\nabla}}
\newcommand{\Div} {\ensuremath{\nabla\cdot}}

\newcommand{\jump}[1] {\ensuremath{\llbracket#1\rrbracket}}
\newcommand{\avg}[1] {\ensuremath{\LRc{\!\{#1\}\!}}}

\newcommand{\Lk}{L^2\LRp{D^k}}

\newcommand{\Gh}{\Gamma_h}

\newcommand{\Oh}{\Omega_h}

\newtheorem{theorem}{Theorem}[section]
\newtheorem{lemma}[theorem]{Lemma}

\newcommand{\eval}[2][\right]{\relax
  \ifx#1\right\relax \left.\fi#2#1\rvert}

\newcommand{\note}[1]{#1}

\newcolumntype{C}[1]{>{\centering\let\newline\\\arraybackslash\hspace{0pt}}m{#1}}

\newcommand*\diff[1]{\mathop{}\!{\mathrm{d}#1}}

\makeatletter
\renewcommand\d[1]{\mspace{6mu}\mathrm{d}#1\@ifnextchar\d{\mspace{-3mu}}{}}
\makeatother

\date{}
\graphicspath{{./figs/}}

\begin{document}

\begin{frontmatter}
\title{On the \note{penalty stabilization mechanism for upwind discontinuous Galerkin formulations of first order hyperbolic systems}}

\author[rice]{Jesse Chan\corref{cor1}}
\ead{Jesse.Chan@caam.rice.edu}
\cortext[cor1]{Principal Corresponding author}
\author[vt]{T. Warburton}
\ead{tcew@vt.edu}
\address[rice]{Department of Computational and Applied Mathematics, Rice University, 6100 Main St, Houston, TX, 77005}
\address[vt]{Department of Mathematics, Virginia Tech, 225 Stanger Street, Blacksburg, VA 24061-0123}

\address{}

\begin{abstract}
Penalty fluxes are dissipative numerical fluxes for high order discontinuous Galerkin (DG) methods which depend on a penalization parameter \cite{warburton2013low, ye2016discontinuous}.  We investigate the dependence of the spectra of high order DG discretizations on this parameter, and show that as its value increases, the \note{spectrum} of the DG discretization splits into two disjoint sets of eigenvalues.  One set converges to the eigenvalues of a conforming discretization, while the other set corresponds to spurious eigenvalues which are damped proportionally to $\tau$.  Numerical experiments also demonstrate that undamped spurious modes present in both in the limit of zero and large penalization parameters are damped for moderate values of the upwind parameter.  
\end{abstract}
\end{frontmatter}

\section{Introduction}

High order discontinuous Galerkin (DG) methods are a popular choice of discretization for simulations of time-domain wave propagation, due to their low numerical dispersion and their ability to accomodate unstructured meshes \cite{grote2006discontinuous,hesthaven2007nodal}.  For PDEs in second order form, a common choice of DG method is the symmetric interior penalty (SIPG) discontinuous Galerkin method \cite{riviere2008discontinuous}, which involves a penalization parameter $\tau$.  The role of this penalty parameter was investigated in \cite{Warburton20063205} for discretizations of the Maxwells equations.  There, it was shown that in the limit of $\tau\rightarrow\infty$, the \note{spectrum} of the DG discretization matrix splits into two disjoint sets of eigenvalues, with one set converging to eigenvalues of a conforming discretization and one set associated with divergent spurious eigenvalues.  

We consider penalty fluxes for discretizations of first order systems of hyperbolic PDEs, which depend on some parameter $\tau \geq 0$.  For $\tau = 0$, these fluxes reduce to the energy-conserving and non-dissipative central fluxes \cite{fezoui2005convergence}.  For $\tau > 0$, penalty fluxes are related to upwind fluxes for problems with continuous coefficients.  We extend the results of \cite{Warburton20063205} to DG formulations for first order systems with penalty fluxes and show that the \note{spectrum} of the DG discretization matrix again splits into two disjoint sets of eigenvalues as $\tau\rightarrow\infty$.  The first set of eigenvalues correspond to the eigenvalues of a conforming discretization, while the second set of eigenvalues have real parts approaching $-\infty$ as $\tau\rightarrow \infty$ and correspond to damped spurious modes \note{which are non-conforming}.  Numerical experiments verify these results for the advection and acoustic wave equations in one and two dimensions.  Numerical experiments also show that, while increasing $\tau$ pushes most eigenvalues further left of the imaginary axis, the real parts of certain eigenvalues converge to zero as $\tau\rightarrow\infty$.  This implies that taking a moderate penalty parameter $\tau = 1$ damps spurious \note{\emph{conforming}} modes which remain otherwise undamped in time-dependent simulations with conforming discretizations.  \note{We show that penalized DG methods address both \textit{non-conforming} spurious modes present when using central fluxes ($\tau = 0$) and \textit{conforming} spurious modes which are present under conforming discretizations (or equivalently, under large penalty parameters $\tau \rightarrow \infty$).}  

\section{Semi-discrete discontinuous Galerkin formulation}

We consider linear first order hyperbolic systems of PDEs in $\mathbb{R}^d$
\begin{equation}
\pd{\bm{U}}{t} + \sum_{i=1}^d \bm{A}_i\pd{{\bm{U}}}{\bm{x}_i} = 0
\label{eq:cons}
\end{equation}
with solution components $\bm{U} = (u_1, \ldots, u_M)^T$ and matrices $\bm{A}_i \in \mathbb{R}^{M\times M}$.  
\note{
For hyperbolic systems, any linear combination of $\bm{A}_i$ 
\begin{equation}
\bm{A}_{\alpha} = \sum_{i=1}^d \bm{\alpha}_i \bm{A}_i
\label{eq:hyperbolic}
\end{equation}
is diagonalizable with real eigenvalues.  
}

We assume the domain $\Omega$ is triangulated exactly using a mesh $\Oh$ consisting of elements $D^k$.  We define the $L^2$ inner product over elements $D^k$ and their boundaries $\partial D^k$ as
\[
\LRp{u,v}_{D^k} = \int_{D^k} u(\bm{x})v(\bm{x}) \diff x, \qquad \LRa{u,v}_{\partial D^k} = \int_{\partial D^k} u(\bm{x})v(\bm{x}) \diff x.
\]
A semi-discrete DG formulation for (\ref{eq:cons}) may be written as 
\begin{equation}
\sum_{D^k \in \Oh} \LRp{\pd{\bm{U}}{t},\bm{V} }_{\Lk} = \sum_{D^k \in \Oh}\LRp{ \note{\LRp{\sum_{i=1}^d\bm{A}_i \bm{U}, \pd{\bm{V}}{\bm{x}_i}}_{\Lk}} - \LRa{(\bm{A}_n\bm{U})^*,\bm{V}}_{\partial D^k}}
\label{eq:formulation}
\end{equation}
where $(\bm{A}_n\bm{U})^*$ is a numerical flux defined on boundary faces and shared faces between elements.  

\note{
In the main proof of this paper, we make the assumption that the matrices $\bm{A}_i$ are symmetric and spatially constant for clarity and simplicity of presentation.  However, the numerical fluxes used are more generally applicable (for example, when $\bm{A}_i$ is non-symmetric, spatially varying, or discontinuous), and we present generalizations beyond symmetric and constant $\bm{A}_i$ where possible.  
}

Finally, we note that
the results in this work are straightforward to generalize to the following modification of (\ref{eq:cons})
\begin{equation}
\bm{A}_0\pd{\bm{U}}{t} + \sum_{i=1}^d \bm{A}_i\pd{{\bm{U}}}{\bm{x}_i} = \bm{f}.
\label{eq:consA}
\end{equation}
We assume $\bm{A}_0 = \bm{I}$ in this work for brevity; however, all results are straightforward to extend to systems of the form for any $\bm{A}_0 = \bm{A}_0(\bm{x})$ which may vary spatially and is pointwise symmetric and positive-definite.  Many models of wave propagation in general heterogeneous media can be represented in the form of (\ref{eq:consA}), including the acoustic wave equation \cite{chan2016weight1}, linear elasticity \cite{de2008interior,wilcox2010high,ye2016discontinuous}, and Maxwell's equations \cite{hesthaven2002nodal,hesthaven2004high,grote2007interior,warburton2013low}.  
\note{Additionally, we note that many hyperbolic systems of interest are symmetrizable, and can be cast into the form of (\ref{eq:consA}) with symmetric $\bm{A}_i$ through an appropriate change of variables \cite{kopriva2014energy}.  }

\subsection{Penalization terms in numerical fluxes}
\label{sec:scaletau}

Let $f$ be a face \note{shared between} an element $D^{k,-}$ and its neighbor $D^{k,+}$.  We denote the outward normal on the face $f$ of $D^{k,-}$ as $\bm{n}$, and let $\bm{U}^+, \bm{U}^-$ be evaluations of $\bm{U}$ restricted to $D^{k,+}$ and $D^{k,-}$, respectively.  The jump and average scalar-valued functions are defined as 
\[
\note{\jump{u} = u^+ - u^-, \qquad \avg{u} = \frac{{u}^+ + u^-}{2}.}
\]
We define the normal matrix across a face as
\[
{\bm{A}}_n = \sum_{i=1}^d {\bm{A}_i\bm{n}_i}.
\]
\note{We first introduce the penalty flux for continuous coefficients (where the value of $\bm{A}_n$ is single valued across shared faces) as}
\[
(\bm{A}_n\bm{U})^* = \bm{A}_n\avg{\bm{U}} - \frac{\tau}{2} \bm{A}_n^T \bm{A}_n\jump{\bm{U}}, \qquad \tau \geq 0.
\]
We compare the penalty flux to the well-known upwind numerical flux.  For \note{continuous coefficients, the upwind flux can be written using an eigenvalue decomposition of $\bm{A}_n$ (see, for example, \cite{ye2016discontinuous} and Chapter 2.4 of \cite{hesthaven2007nodal})}.  
Under assumption (\ref{eq:hyperbolic}), ${\bm{A}}_n$ has real eigenvalues, and admits an eigenvalue decomposition
\[
\bm{A}_n = \bm{V}{\bm{\Lambda}}\note{\bm{V}^{-1}}, \qquad \bm{\Lambda} = 
\left(\begin{array}{ccc}
\lambda_1 & & \\
& \ddots & \\
& & \lambda_d
\end{array}\right)
\]
For problems with continuous coefficients, the upwind numerical flux over a face $f \in \Gh$ can then be defined  as
\[
(\bm{A}_n\bm{U})^* = \bm{A}_n^+\bm{U}^- + \bm{A_n}^- \bm{U}^+,
\]
where the matrices $\bm{A_n}^+,\bm{A_n}^-$ are constructed from the positive and negative eigenvalues 
\begin{align*}
\bm{A}_n^+ &= \frac{1}{2}\bm{V} \LRp{\bm{\Lambda} + \LRb{\bm{\Lambda}}} \note{\bm{V}^{-1}}\\
\bm{A}_n^- &= \frac{1}{2}\bm{V} \LRp{\bm{\Lambda} - \LRb{\bm{\Lambda}}} \note{\bm{V}^{-1}},
\end{align*}
and $\LRb{\bm{\Lambda}}$ is the diagonal matrix whose entries consist of the absolute values of the eigenvalues $\LRb{\lambda_i}$.  This can be rewritten as the sum of the central flux and a penalization
\[
(\bm{A}_n\bm{U})^* = \bm{A}_n\avg{\bm{U}} - \frac{1}{2}\bm{V}\LRb{\bm{\Lambda}}\note{\bm{V}^{-1}} \jump{\bm{U}}.  
\]
\note{
Assuming constant, symmetric $\bm{A}_i$ and a penalty numerical flux, one can show the semi-discrete energy estimate \cite{hesthaven2007nodal, warburton2013low, kopriva2014energy, chan2016weight1} for (\ref{eq:cons})
\[
\sum_{D^k\in \Oh}\pd{}{t}\nor{\bm{U}}^2_{L^2\LRp{D^k}} = -\sum_{D^k \in \Oh} \frac{\tau}{2}\nor{\bm{A}_n\jump{\bm{U}}}^2_{L^2\LRp{\partial D^k}} \leq 0.
\]
The penalization term can thus be interpreted as adding a stabilizing contribution which dissipates components of $\bm{U}$ for which $\bm{A}_n\jump{\bm{U}}$ is large, with larger values of $\tau$ increasing the rate of dissipation (this is true even when $\bm{A}_n$ is not symmetric).  We will relate this dissipation to the spectrum of the DG discretization matrix in Section~\ref{sec:spec}.  In particular, we show that the real parts of eigenvalues corresponding to components of eigenvectors with $\bm{A}_n\jump{\bm{U}}\neq 0$ approach $-\infty$ as $\tau \rightarrow \infty$.  

When $\bm{A}_n$ is not continuous across faces, the penalty numerical flux can be generalized as follows \cite{ye2016discontinuous}
\[
(\bm{A}_n\bm{U})^* = \avg{\bm{A}_n\bm{U}} - \frac{\tau}{2} \bm{A}_n^T \jump{\bm{A}_n\bm{U}},
\]
where we define the normal jump and average
\[
\jump{\bm{A}_n\bm{U}}\coloneqq \sum_{i=1}^d\bm{n}_i \jump{ \bm{A}_i \bm{U}}, \qquad \avg{\bm{A}_n\bm{U}} \coloneqq \sum_{i=1}^d\bm{n}_i\avg{ \bm{A}_i \bm{U}}.
\]
If $\bm{A}_i$ is constant over each element $D^k$, we again have an energy estimate with a dissipative contribution from the penalty flux
\[
\sum_{D^k\in \Oh}\pd{}{t}\nor{\bm{U}}^2_{L^2\LRp{D^k}} = -\sum_{D^k \in \Oh} \frac{\tau}{2}\nor{\jump{\bm{A}_n\bm{U}}}^2_{L^2\LRp{\partial D^k}} \leq 0.
\]
Finally, we note that the energy estimate changes when $\bm{A}_i$ varies over an element, though the penalty flux still adds a dissipative contribution to the right hand side \cite{kopriva2014energy}.  
}

\subsection{Scaling of $\tau$}

We note that including a positive-definite weight matrix $\bm{W}$ into the penalty term
\[
\frac{\tau}{2} \bm{A}_n^T \bm{W}\bm{A}_n\jump{\bm{U}}
\]
recovers the upwind fluxes for $\tau = 1$ and $\bm{W} = \bm{V}\LRb{\bm{\Lambda}}^{\dagger}\bm{V}^{-1}$, where $\bm{A}^{\dagger}$ is the pseudo-inverse of $\bm{A}$ \cite{horn2012matrix}.  However, the advantage of penalty fluxes is that, unlike upwind fluxes, they can be used when explicit knowledge of a \note{joint} diagonalization of $\bm{A}_n$ is not known, such as for linear elasticity in anisotropic media \cite{ye2016discontinuous}.  Thus, we assume $\bm{W} = \bm{I}$ in this work.  

\note{It remains to specify an appropriate value of $\tau$.  We first note that the timestep is proportional to the spectral radius of the DG matrix, which increases proportionally to $\tau$ \cite{chan2015gpu}.  Since both upwind and penalty fluxes are identical apart from the penalization term, we choose $\tau$ such that the penalty flux and upwind penalization terms are of the same magnitude 
\[
\nor{\frac{\tau}{2} \bm{A}_n^T\bm{A}_n} = \nor{\frac{1}{2}\bm{V}\LRb{\bm{\Lambda}}{\bm{V}^{-1}}}.
\]
This is sufficient to ensure that the spectral radius of the DG matrix is similar when using either upwind or penalty fluxes.  

The penalty flux and upwind penalization terms can both be bounded from above in terms of the largest magnitude eigenvalue $\max_i\LRb{\lambda_i}$ and the condition number $\kappa\LRp{\bm{V}}$ of the eigenvector matrix $\bm{V}$
\[
\nor{\frac{1}{2}\bm{V}\LRb{\bm{\Lambda}}{\bm{V}^{-1}}} \leq \max_i \LRb{\lambda_i} \kappa\LRp{\bm{V}}, \qquad \nor{\bm{A}_n^T\bm{A}_n} = \nor{\bm{V}^{-T}\bm{\Lambda}\bm{V}^T \bm{V}\bm{\Lambda}\bm{V}^{-1}} \leq \max_i\LRb{\lambda_i}^2\kappa\LRp{\bm{V}}^2.
\]
Thus, choosing $\tau = 1/\LRp{\max_i\LRb{\lambda_i}\kappa\LRp{\bm{V}}}$ results in a penalization term whose magnitude is the same as that of the upwind penalization term.  When $\bm{A}_n$ is symmetric (or symmetrizable through a change of variables), $\bm{V}$ is unitary and $\kappa\LRp{\bm{V}} = 1$, simplifying the selection of $\tau$ to $\tau = 1/\max_i\LRb{\lambda_i}$.  
}

\section{Dependence of \note{spectrum} on the penalty parameter}
\label{sec:spec}
The remainder of this work focuses on the influence of $\tau$ on the \note{spectrum} of DG discretization matrices.  We will show that, as $\tau$ increases, the \note{spectrum} splits into two sets of eigenvalues corresponding respectively to conforming and spurious modes, where conformity is defined based on $\bm{A}_n$.  \note{For clarity and simplicity of presentation, we will assume that $\bm{A}_i$ are symmetric and spatially constant for the remainder of this paper.  }

\subsection{Conforming and non-conforming approximation spaces}

Let $V$ denote the piecewise polynomial approximation space \note{over an element $D^k$}
\[
V = \LRc{ \bm{U} \in\LRp{ L^2(\Omega)}^d : \left.\bm{U}\right|_{D^k} \in \LRp{P^N(D^k)}^d, \quad \note{\forall D^k \in \Oh}}.  
\]
We decompose $V$ into a conforming space $V^C$ and a non-conforming space $V^{NC}$ based on the penalization term.  The conforming space is defined to be elements of $V$ for which the penalization term vanishes over all faces of all elements.  For penalty fluxes, \note{this term is $\bm{A}_n^T \bm{A}_n\jump{\bm{U}}$, such that} $V^C$ is
\begin{equation}
V^{C} = \LRc{ \bm{U}(\bm{x}) \in V: \quad  \bm{A}_n\jump{\bm{U}} = 0, \quad \forall f \in \partial D^k, \quad \forall D^k\in \Oh}.
\label{eq:conf}
\end{equation}
\note{This simplification can be made due to the fact that the null space of $\bm{A}_n^T\bm{A}_n$ is identical to the null space of $\bm{A}_n$.\footnote{\note{The null space of $\bm{A}_n^T\bm{A}_n$ contains the null space of $\bm{A}_n$.  However, if $\bm{A}_n^T\bm{A}_n\bm{U} = 0$, then $\bm{U}^T\bm{A}_n^T\bm{A}_n\bm{U} = \nor{\bm{A}_n\bm{U}}^2 = 0$ as well, implying that the null space of $\bm{A}_n$ also contains the null space of $\bm{A}_n^T\bm{A}_n$.
}}
}  
The non-conforming approximation space $V^{NC}$ is defined as $L^2$ orthogonal complement of $V^{C}$ in $V$
\[
V^{NC} = \LRc{ \bm{U}(\bm{x}) \in V: \LRp{\bm{u},\bm{v}}_{L^2(\Omega)} = 0, \quad \forall \bm{v}\in V^C}.
\]
A conforming space could also be defined using the upwind flux by seeking $\bm{U}$ such that $\bm{V}\LRb{\bm{\Lambda}}\note{\bm{V}^{-1}} \jump{\bm{U}} = 0$.  When the coefficient matrices $\bm{A}_i$ are spatially continuous, the condition $\bm{V}\LRb{\bm{\Lambda}}\note{\bm{V}^{-1}} \jump{\bm{U}} = 0$ and the penalty flux condition $\bm{A}_n^T \bm{A}_n\jump{\bm{U}} = 0$ are both equivalent to
\[
\bm{v}_i^T \jump{\bm{U}} = 0, \qquad \lambda_i \neq 0,
\]
where $\lambda_i, \bm{v}_i$ are an eigenvalue and eigenvector of $\bm{A}_n$.  As a result, the conforming and non-conforming spaces induced by upwind and penalty fluxes are identical.  

We note that, for problems with discontinuous material data, it is less straightforward to show equivalence between the conforming spaces $V^C$ resulting from upwind and penalty fluxes.  However, penalty fluxes can still be applied in the presence of discontinuous coefficients by incorporating these spatial variations into the matrix $\bm{A_0}$.  When paired with penalty fluxes, this approach results in energy stable DG methods \cite{chan2016weight1}.  Additionally, DG methods with penalty fluxes are observed to perform similarly to DG methods with upwind fluxes for several realistic wave propagation problems with discontinuous coefficients \cite{ye2016discontinuous}.  

\subsection{Scalar advection and acoustic wave propagation}
\label{sec:confexamples}
We illustrate the conforming spaces induced by penalty fluxes using the scalar advection equation and acoustic wave equation as examples.  The scalar advection equation is given as
\[
\pd{u}{t} + \Div\LRp{\bm{\beta}u} = 0.
\]
where $\bm{\beta}$ is the direction of advection.  \note{We also assume $\bm{\beta}$ is constant, such that the resulting equation can be rearranged into the form (\ref{eq:cons}).}  The penalty flux is 
\[
(\bm{A}_n\bm{U})^* = \bm{\beta}_n\avg{ u} - \tau\frac{\LRp{\bm{\beta}_n}^2}{2}\jump{u}.
\]
\note{Setting $\tau = \LRb{\bm{\beta}_n}$ recovers the standard upwind flux.}  
The conforming space $V^C$ induced by this flux is then
\[
V^C = \LRc{ u \in L^2\LRp{\Omega} : \quad \left.u\right|_{D^k} \in P^N(D^k), \quad \LRb{\bm{\beta}_n}\jump{u} = 0}.
\]
In one dimension, this simply implies that $u$ is $C^0$-continuous across element boundaries.  In higher dimensions, this implies that $u$ is continuous along streamlines or directions $\bm{d}(\bm{x})$ where $\bm{\beta}\cdot \bm{d} \neq 0$.  

Next, we consider the acoustic wave equation in pressure-velocity form
\begin{align*}
\frac{1}{c^2}\pd{p}{t} + \Div \bm{u} &= 0\\
\pd{\bm{u}}{t} + \Grad p &= 0,
\end{align*}
where $c^2(\bm{x})$ is the wavespeed.  Let $\bm{U}$ denote the group variable $\bm{U} = (p,u,v)$, where $u$ and $v$ are the $x$ and $y$ components of velocity.  Then, in two dimensions, the isotropic wave equation is given as
\[
\pd{\bm{U}}{t} + \pd{\bm{A}_x\bm{U}}{x} + \pd{\bm{A}_y\bm{U}}{y} = 0, \qquad \bm{A}_x = 
\left(\begin{array}{ccc}
0 & 1 & 0\\
1 & 0 & 0\\
0 & 0 & 0
\end{array}
\right), \qquad 
\bm{A}_y = 
\left(\begin{array}{ccc}
0 & 0 & 1\\
0 & 0 & 0\\
1 & 0 & 0
\end{array}
\right).
\]
The normal flux matrix $\bm{A}_n$ is then
\[
\bm{A}_n = 
\left(\begin{array}{ccc}
0 & \bm{n}_x & \bm{n}_y\\
\bm{n}_x & 0 & 0\\
\bm{n}_y & 0 & 0
\end{array}
\right)
\]
implying that the penalty fluxes are 
\[
\bm{A}_n \avg{\bm{U}} - \frac{\tau}{2}\bm{A}_n^T\bm{A}_n \jump{\bm{U}} = \left(\begin{array}{c}
\avg{\bm{u}_n}\\
\avg{p }\bm{n}_x\\
\avg{p }\bm{n}_y
\end{array}
\right) - 
\frac{\tau}{2}\left(\begin{array}{c}
\jump{p}\\
\jump{\bm{u_n}}\bm{n}_x\\
\jump{\bm{u_n}}\bm{n}_y
\end{array}
\right)
\]
Since the eigenvalues of $\bm{A}_n$ are $-1,0,1$, we take the penalty parameter to be $\tau = 1 / (\max_i |\lambda_i|) = 1$.   Additionally, for $\tau = 1$, the penalty fluxes coincide with the upwind fluxes for continuously varying media.  In this case, the conforming subspace $V^{C}$ induced by the penalty flux consists of a scalar (pressure) component $V^C_p$ and a vector (velocity) component $V^C_{\bm{u}}$.  The pressure component $p$ satisfies $\jump{p} = 0$ over all faces $f\in \Gh$.  For polynomial approximation spaces, this implies that $p$ is continuous across faces, edges, and vertices, and $V^C_p$ is the standard $C^0$ continuous piecewise polynomial finite element space.  For the velocity component $V_{\bm{u}}^C$, the penalty fluxes enforce normal continuity, such that $\jump{\bm{u}}\cdot \bm{n} = 0$ over all faces $f \in \Gh$.  Since each component of $\bm{u}$ is approximated from $P^N(D^k)$, this implies that $V^C_{\bm{u}}$ is the $H({\rm div})$-conforming Brezzi-Douglas-Marini finite element space \cite{brezzi1985two,boffi2013mixed} \note{on triangular and tetrahedral elements}.\footnote{It is possible to recover $V^C_{\bm{u}}$ as the Raviart-Thomas finite element space by approximating each component of $\bm{u}$ from a different polynomial space \cite{kirby2004algorithm}.}   These reflect minimal continuity conditions which result from requiring $\Div \bm{u} \in \LRp{L^2\LRp{\Oh}}^d$.

Penalty fluxes can also be interpreted as Lax-Friedrichs fluxes applied to scalar and vector variables separately.  In contrast, for a component-wise Lax-Friedrichs flux 
\[
(\bm{A}_n\bm{U})^* = \left(\begin{array}{c}
\avg{\bm{u}_n}\\
\avg{p }\bm{n}_x\\
\avg{p }\bm{n}_y
\end{array}
\right) - 
\frac{\tau}{2}\left(\begin{array}{c}
\jump{p}\\
\jump{\bm{u}_x}\\
\jump{\bm{u}_y}
\end{array}
\right),
\]
the conforming velocity space $V^C_{\bm{u}}$ becomes the space of vector piecewise polynomial functions which are $C^0$ continuous in both $x$ and $y$ components, resulting in over-constrained continuity conditions for $\bm{u}$.

\subsection{Behavior as $\tau\rightarrow\infty$}

Let $\bm{K}$ be the matrix resulting from the discretization of the DG formulation (\ref{eq:formulation}).  We are interested in eigenvalues and eigenvectors of 
\[
\bm{K}\bm{u} = \lambda_i\bm{M}\bm{u}
\]
where $\bm{M}$ is the $L^2$ mass matrix.  We adapt the approach of Warburton and Embree in \cite{Warburton20063205} to (\ref{eq:block}) to show that the \note{spectrum} of $\bm{K}$ splits into two sets of eigenvalues as $\tau\rightarrow \infty$, the first of which approaches the eigenvalues of $\bm{A}$, and the second of which diverges with real part approaching $-\infty$ as $\tau\rightarrow \infty$.  

\note{We first note that integrating by parts the right hand side of (\ref{eq:formulation}) and manipulating interface terms yields the following rearrangement of the DG formulation
\begin{align}
b(\bm{U},\bm{V}) &=\sum_{D^k\in \Oh} \frac{1}{2}\LRp{ \LRp{\sum_{i=1}^d\bm{A}_i\bm{U}, \pd{\bm{V}}{\bm{x}_i}}_{D^k} - \LRp{\sum_{i=1}^d\bm{A}_i \pd{\bm{U}}{\bm{x}_i}, {\bm{V}}}_{D^k}} \label{eq:dgformskew}\\
&+ \sum_{D^k\in \Oh}\frac{1}{2}\LRp{\LRa{\bm{A}_n\jump{\bm{U}}, \avg{\bm{V}}}_{\partial D^k}-\LRa{\avg{\bm{U}}, \bm{A}_n\jump{\bm{V}}}_{\partial D^k}} - \frac{\tau}{2}\LRa{\bm{A}_n\jump{\bm{U}},\bm{A}_n\jump{\bm{V}}}_{\partial D^k}.
\nonumber
\end{align}
The formulation is skew-symmetric (assuming $\bm{A}_i$ are constant and symmetric), with the exception of the symmetric semi-definite $\tau$-dependent penalization term.}

Let $\bm{K}$ denote the DG discretization matrix resulting from (\ref{eq:dgformskew}).  Let $\LRc{\bm{\Phi}_i }_{i=1}^{N^{C}}$ and $\LRc{\bm{\Psi}_i }_{i=1}^{N^{NC}}$ now be $L^2$-orthogonal bases for $\bm{V}^C$ and $\bm{V}^{NC}$, respectively.  These spaces induce a block decomposition of $\bm{K}$
\begin{equation}
\bm{K} = \left(\begin{array}{cc}
\bm{A} & \bm{B}\\
-{\bm{B}}^T & \bm{C} + \tau \bm{S}
\end{array}\right),
\label{eq:block}
\end{equation}
\note{
where the blocks $\bm{A},\bm{B},\bm{C}$ are defined as
\begin{align*}
\bm{A}_{mn} &= \frac{1}{2}\sum_{D^k\in \Oh} \LRp{ \sum_{i=1}^d\LRp{\LRp{\bm{A}_i\bm{\Phi}_n, \pd{\bm{\Phi}_m}{\bm{x}_i}}_{D^k} - \LRp{\bm{A}_i \pd{\bm{\Phi}_n}{\bm{x}_i}, {\bm{\Phi}_m}}_{D^k}}}, \qquad 1\leq m,n \leq N^C\\
\bm{B}_{mn} &= \frac{1}{2}\sum_{D^k\in \Oh} \LRp{ \sum_{i=1}^d\LRp{\LRp{\bm{A}_i\bm{\Phi}_n, \pd{\bm{\Psi}_m}{\bm{x}_i}}_{D^k} - \LRp{\bm{A}_i \pd{\bm{\Phi}_n}{\bm{x}_i}, {\bm{\Psi}_m}}_{D^k}} - {\LRa{\avg{\bm{\Phi}_n}, \bm{A}_n\jump{\bm{\Psi}_m}}_{\partial D^k}} }, \quad 1\leq n \leq N^{NC} \\
\bm{C}_{mn}  &=  \frac{1}{2}\sum_{D^k\in \Oh}\LRp{
\sum_{i=1}^d\LRp{\LRp{\bm{A}_i\bm{\Psi}_n, \pd{\bm{\Psi}_m}{\bm{x}_i}}_{D^k} - \LRp{\bm{A}_i \pd{\bm{\Psi}_n}{\bm{x}_i}, {\bm{\Psi}_m}}_{D^k}}}\\
&\qquad + \frac{1}{2}\sum_{D^k\in \Oh}\LRp{\LRa{\bm{A}_n\jump{\bm{\Psi}_n}, \avg{\bm{\Psi}_m}}_{\partial D^k}-\LRa{\avg{\bm{\Psi}_n}, \bm{A}_n\jump{\bm{\Psi}_m}}_{\partial D^k}}, \qquad 1 \leq m,n \leq N^{NC}\\
\bm{S}_{mn}  &= -\sum_{D^k\in \Oh}  \frac{1}{2}\LRa{\bm{A}_n\jump{\bm{\Psi}_n},\bm{A}_n\jump{\bm{\Psi}_m}}_{\partial D^k}, \qquad 1 \leq m,n \leq N^{NC}.  
\end{align*}
where we have used that, because $\bm{\Phi}_n\in V^{C}$, $\bm{A}_n\jump{\bm{\Phi}_n} = 0$.}

\note{We note that $\bm{A}$ is skew-symmetric, and $\bm{S}$ is symmetric and negative semi-definite.  To show that $\bm{S}$ is also negative definite, assume that there exists some $\bm{v} \in V^{NC}$ with degrees of freedom $\bm{V}$ such that $\bm{S}\bm{V} = 0$ and
\[
\bm{V}^T\bm{S}\bm{V} = -\frac{1}{4}\sum_{D^k \in \Oh} \int_{\partial D^k} \LRp{\bm{A}_n\jump{ \bm{v}_n}}^2 = 0.
\]
However, because each term in the summation on the right hand side is negative definite, we conclude that $\bm{A}_n\jump{ \bm{v}_n} = 0$ across all mesh interfaces.  This would imply that $\bm{v}\in V^C$, contradicting the assumption that $\bm{v} \in V^{NC}$.}  

We may now show how the \note{spectrum} of $\bm{K}$ behaves as $\tau \rightarrow \infty$.  Since $\bm{A}$ is skew-symmetric and $\bm{S}$ is symmetric, they are diagonalizable under unitary matrices $\bm{U}$ and $\bm{Q}$, whose columns contain the eigenvectors of $\bm{A}$ and $\bm{S}$, respectively.   
Following \cite{Warburton20063205}, we apply a block diagonal similarity transform to yield
\[
\tilde{\bm{K}} = \left(\begin{array}{cc}
\bm{U}^* & \\
& \bm{Q}^*
\end{array}\right)
\left(\begin{array}{cc}
\bm{A} & \bm{B}\\
-\bm{B}^T & \bm{C} + \tau \bm{S}
\end{array}\right)
\left(\begin{array}{cc}
\bm{U} & \\
& \bm{Q}
\end{array}\right)
 = \left(\begin{array}{cc}
\bm{\Lambda}^C & \bm{U}^*\bm{B}\bm{Q}\\
-\bm{Q}^*\bm{B}^T\bm{U} & \bm{Q}^*\bm{C}\bm{Q} + \tau \bm{\Lambda}^{S}
\end{array}\right),
\]
where $\bm{\Lambda}^C,\bm{\Lambda}^{S}$ are diagonal matrices whose entries consist of eigenvalues of $\bm{A}$ and $\bm{S}$.  

Since $\bm{Q}^*\bm{C}\bm{Q}$ is skew-symmetric and its diagonal has zero real part, $\bm{Q}^*\bm{C}\bm{Q} + \tau \bm{\Lambda}^{S}$ has diagonal entries $\tau \lambda^S_j + \mi\gamma_j$, where $\lambda^S_j$ is real and negative and $\gamma_j$ is independent of $\tau$.  Thus, the diagonal of $\tilde{\bm{K}}$ consists of $N^C$ purely imaginary values and $N^{NC}$ values of the form $\tau \lambda^S_j + \mi\gamma_j$.  Since $\bm{B}$ is independent of $\tau$, the entries of $\bm{U}^*\bm{B}\bm{Q}$ are independent of $\tau$, assuming that $\bm{U}$ and $\bm{Q}$ are normalized.  Gerschgorin's theorem applied to $\tilde{\bm{K}}$ then implies that the eigenvalues of $\bm{K}$ are contained in two sets of discs with radii independent of $\tau$.  The first set of discs are centered around $\lambda^C_j$ for $j = 1,\ldots,N^C$, while the second set of discs are centered around $\tau\lambda^S_j + \mi\gamma_j$ for $j = 1,\ldots,N^{NC}$. For sufficiently large $\tau$, these disks must become disjoint, implying that the union of disks centered around  $\lambda^C_i$ contains exactly $N^C$ eigenvalues, while the union of disks centered around  $\lambda^{NC}_i$ contain exactly $N^{NC}$ eigenvalues \cite{horn2012matrix}.  In other words, the \note{spectrum} of $\bm{K}$ diverges into two distinct sets of $N^C$ and $N^{NC}$ eigenvalues each as $\tau \rightarrow \infty$, with the set of $N^{NC}$ eigenvalues containing negative real parts with magnitude $O(\tau)$.  

Consider now the $N^C$ eigenvalues \note{ $\lambda_1, \ldots, \lambda_{N^C}$} of smallest magnitude.  By the Gerschgorin argument, these must remain bounded as $\tau \rightarrow \infty$.  \note{Let $\bm{W}$ denote the matrix of eigenvectors and $\bm{\Lambda}^C$ the matrix of eigenvalues corresponding to these $N^C$ smallest eigenvalues, such that
\[
\bm{W} = \LRp{\begin{array}{c}
\bm{W}^C\\
\bm{W}^{NC}
\end{array}}, \qquad 
\bm{\Lambda}^C = \LRp{\begin{array}{cccc}
\lambda_1 & & &\\
& \lambda_2 & &\\
& & \ddots &\\
& & & \lambda_{N^C}
\end{array}
}.
 \]
Here, $\bm{W}^C, \bm{W}^{NC}$ refer to components of eigenvectors from the conforming and non-conforming spaces $V^C, V^{NC}$, respectively.  These eigenvector matrices have $N^C$ columns, and $\bm{W}^C, \bm{W}^{NC}$ have $N^{C}$ and $N^{NC}$ rows, respectively.  By definition, $\bm{W}^C, \bm{W}^{NC}$ satisfy}
\[
\left(\begin{array}{cc}
\bm{A} & \bm{B}\\
-\bm{B}^T & \bm{C} + \tau \bm{S}
\end{array}\right)
\left(\begin{array}{c}
\bm{W}^C
\\
\bm{W}^{NC}
\end{array}\right) = 
\left(\begin{array}{c}
\bm{A}\bm{W}^C + \bm{B}\bm{W}^{NC}\\
-\bm{B}^T\bm{W}^C + \bm{C}\bm{W}^{NC} + \tau \bm{S}\bm{W}^{NC}
\end{array}\right)
= 
\left(\begin{array}{c}
\bm{W}^C
\\
\bm{W}^{NC}
\end{array}\right) 
\note{\bm{\Lambda}_C}.
\]
This implies 
\[
\tau \nor{\bm{S}\bm{W}^{NC}} = \nor{\bm{\Lambda}^C \bm{W}^{NC} + \bm{B}^T\bm{W}^C  - \bm{C}\bm{W}^{NC}}.
\]
Since $\bm{\Lambda}^C$ remains bounded as $\tau\rightarrow \infty$, the right hand side is bounded independently of $\tau$ \note{(assuming $\bm{W}$ is normalized)}, implying that the non-conforming component $\bm{W}^{NC}$ satisfies
\[
\sqrt{\lambda^{NC}_{\min}}\nor{\bm{W}^{NC}} \leq \nor{\bm{S}\bm{W}^{NC}} = O(1/\tau)
\]
where $\lambda^{NC}_{\min} > 0$ is the smallest eigenvalue of $\bm{S}$.  As a consequence, $\bm{W}^{NC} \rightarrow 0$ as $\tau\rightarrow \infty$, and the smallest $N^C$ eigenvalues of $\bm{K}$ converge to the eigenvalues of $\bm{A}$ at a rate of $O(1/\tau)$.  These eigenvalues of $\bm{A}$ correspond to a discretization using the conforming approximation space (\ref{eq:conf}). 

These results are summarized in the following lemma:
\begin{lemma}
As $\tau \rightarrow \infty$, the spectrum of the DG discretization matrix $\bm{K}$ decouples into two sets of eigenvalues $\LRc{\lambda^C_1,\ldots,\lambda^C_{N^C}}$ and $\LRc{\lambda^{NC}_1,\ldots,\lambda^{NC}_{N^{NC}}}$.  The eigenvalues $\LRc{\lambda^{NC}_1,\ldots,\lambda^{NC}_{N^{NC}}}$ diverge towards the left half plane, with $\LRb{{\rm Re}(\lambda_i)} = O(\tau)$, while  the eigenvalues $\LRc{\lambda^C_1,\ldots,\lambda^C_{N^C}}$ converge to the eigenvalues of $\bm{A}$ with rate $O(1/\tau)$.  
\label{lemma:eig}
\end{lemma}
We note that, while we consider skew-symmetric DG formulations in this work, Lemma~\ref{lemma:eig} is straightforward to generalize to any DG formulations where $\bm{S}$ is symmetric negative definite and $\bm{A}$ is diagonalizable.

\section{Numerical experiments}

Lemma~\ref{lemma:eig} illustrates the behavior of the \note{spectrum} of the DG discretization for the asymptotic cases when $\tau = 0$ or $\tau \rightarrow \infty$.  However, it is less clear how the \note{spectrum} of $\bm{K}$ behaves for $\tau \approx 1$.  We rely instead on numerical experiments in one and two space dimensions to illustrate behaviors for the advection and acoustic wave equations.  {\note{Since $\tau$ can be chosen to be any arbitrary non-negative value, we consider three cases: $\tau = 0$, $\tau = 1$, and $\tau \rightarrow \infty$.  As mentioned in Section~\ref{sec:scaletau}, it may be beneficial to scale $\tau$ with problem parameters.  However, for the simple model problems considered in this section, taking $\tau = 1$ recovers an upwind flux, and the resulting spectrum of the DG operator is clearly distinct from the cases when $\tau = 0$ or $\tau \rightarrow \infty$.}  }

These experiments indicate that certain modes with negative real part for $\tau = 1$ return to the imaginary axis as $\tau \rightarrow \infty$.  Since eigenvalues with a negative real part correspond to dissipation in time-domain simulations, this implies that these under-resolved modes become undamped as $\tau$ grows sufficiently large.

\subsection{1D experiments}

\begin{figure}
\centering
\subfloat[Paths of eigenvalues as $\tau$ increases]{\includegraphics[width=.48\textwidth]{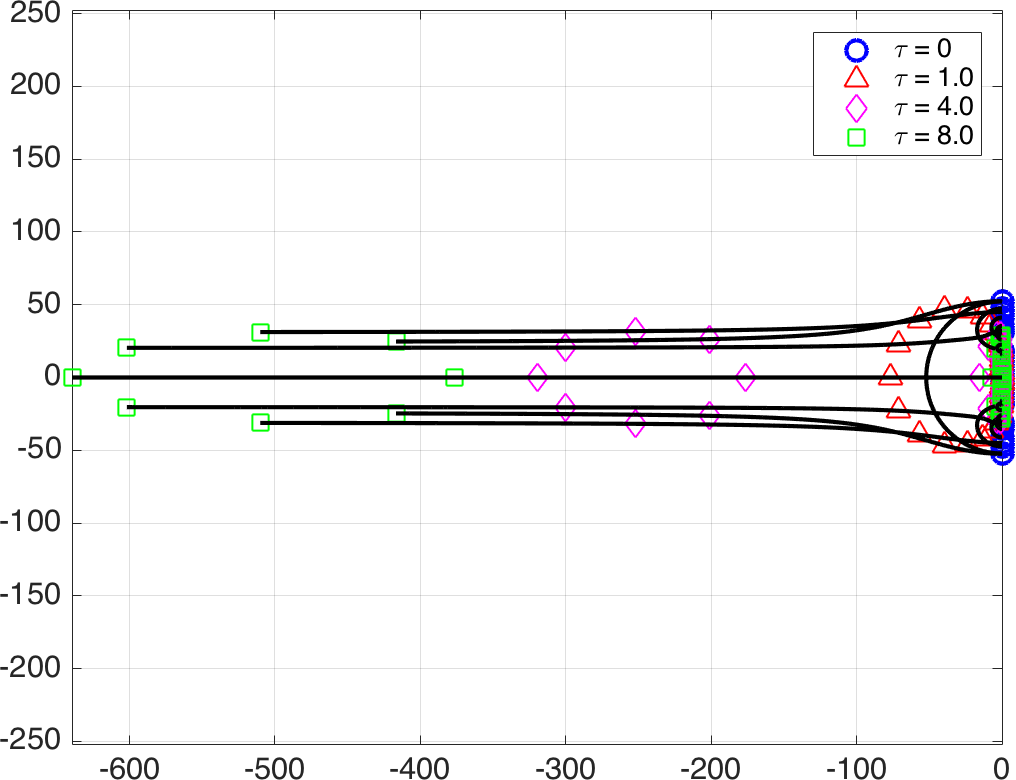}}
\hspace{1em}
\subfloat[Zoom of eigenvalues near imaginary axis]{\includegraphics[width=.435\textwidth]{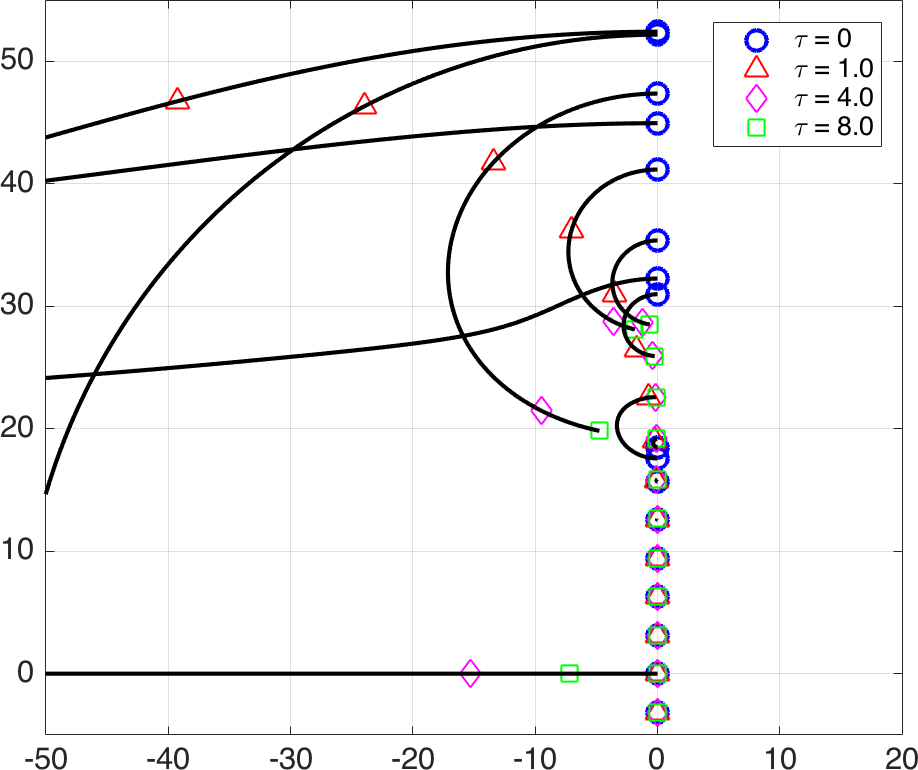}}
\caption{Eigenvalue paths for DG advection with $\tau \in [0,4]$ on a mesh of 8 elements of degree $N=3$.  Eigenvalues are overlaid on these paths for $\tau = 0$, $\tau = 1$, and $\tau = 4$.  The zoomed in view near the imaginary axis shows the return of spurious modes to the imaginary axis for $\tau$ sufficiently large. }
\label{fig:track1D}
\end{figure}

\begin{figure}
\centering
\subfloat{\includegraphics[width=.32\textwidth]{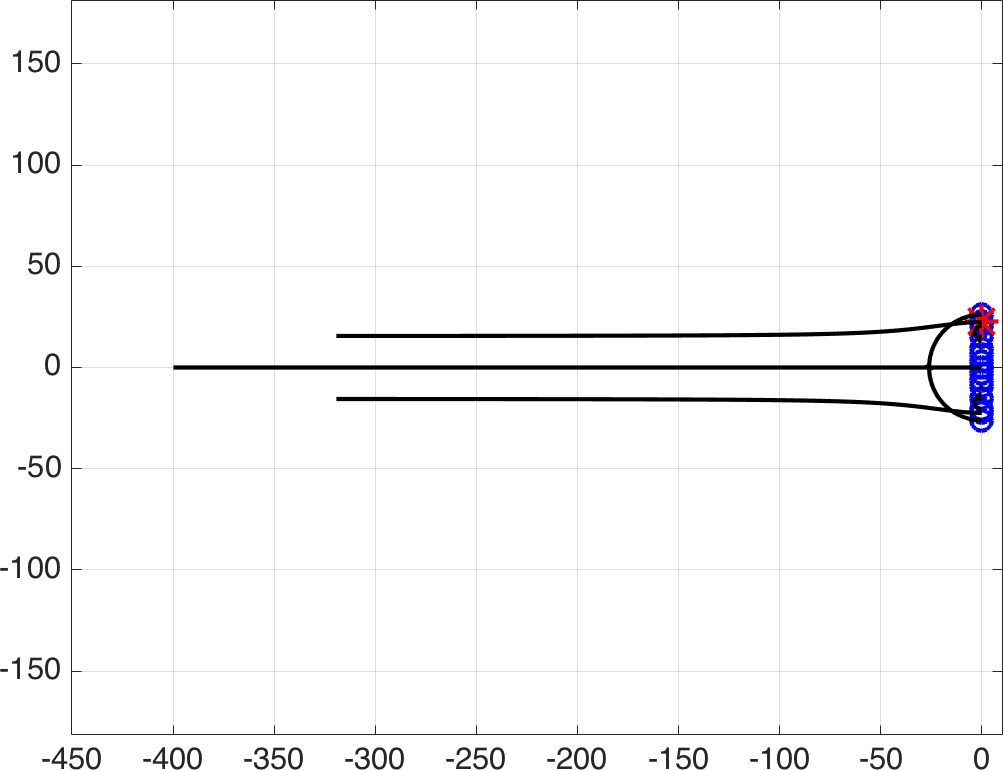}}
\hspace{.5em}
\subfloat{\includegraphics[width=.32\textwidth]{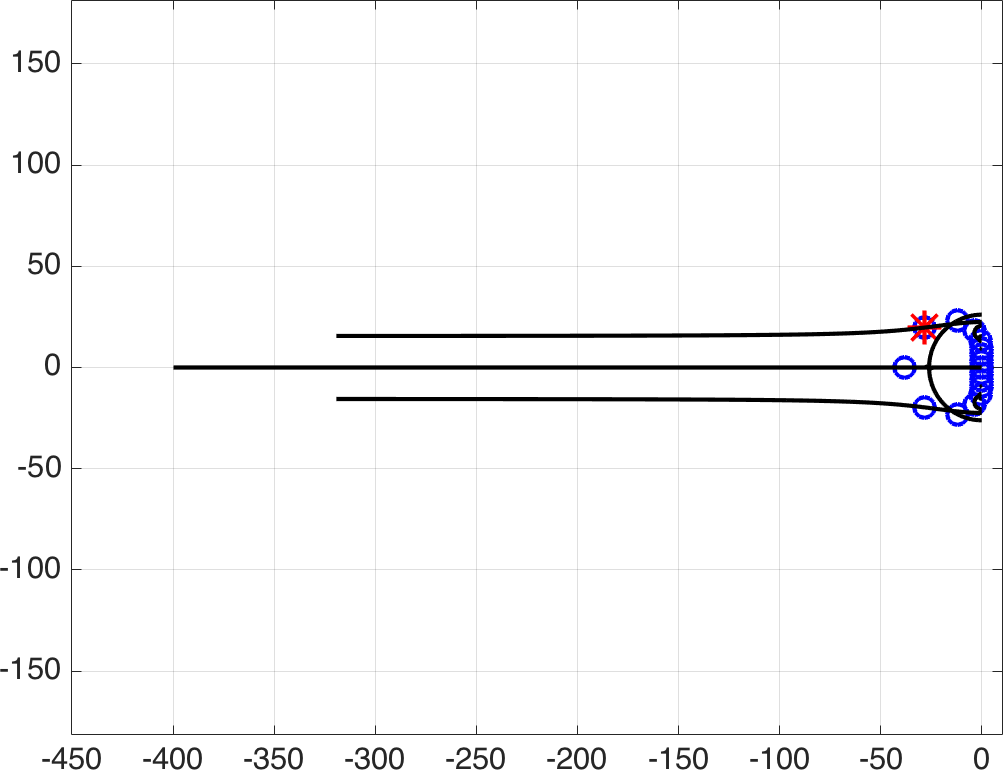}}
\hspace{.5em}
\subfloat{\includegraphics[width=.32\textwidth]{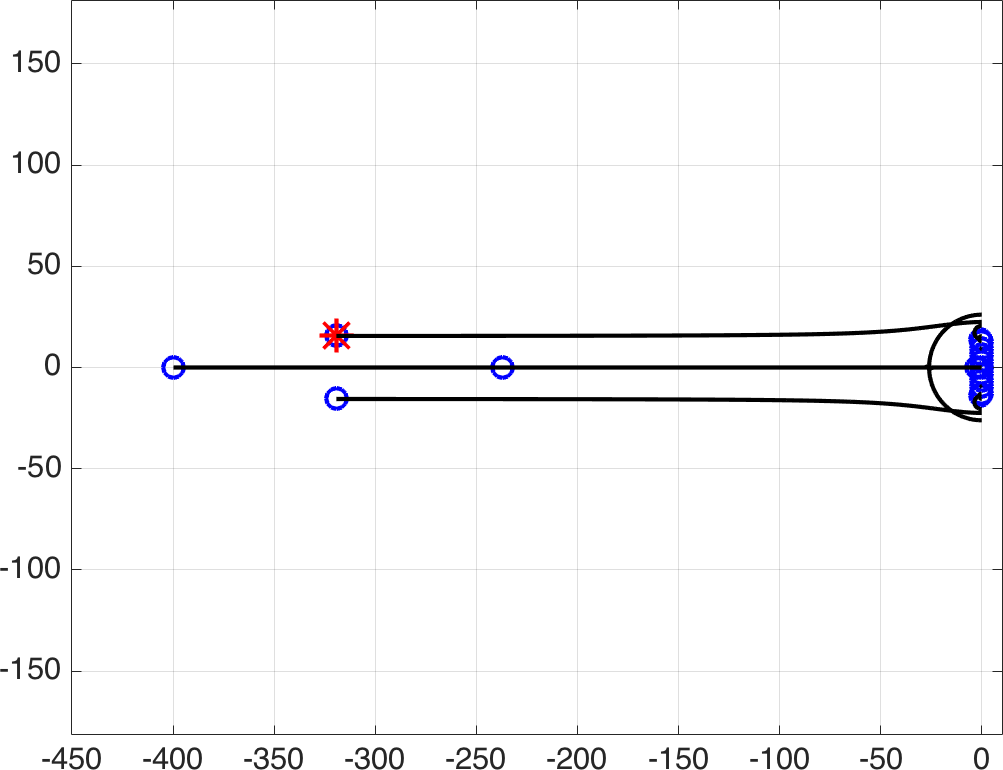}}\\
\subfloat[$\tau = 0$]{\includegraphics[width=.32\textwidth]{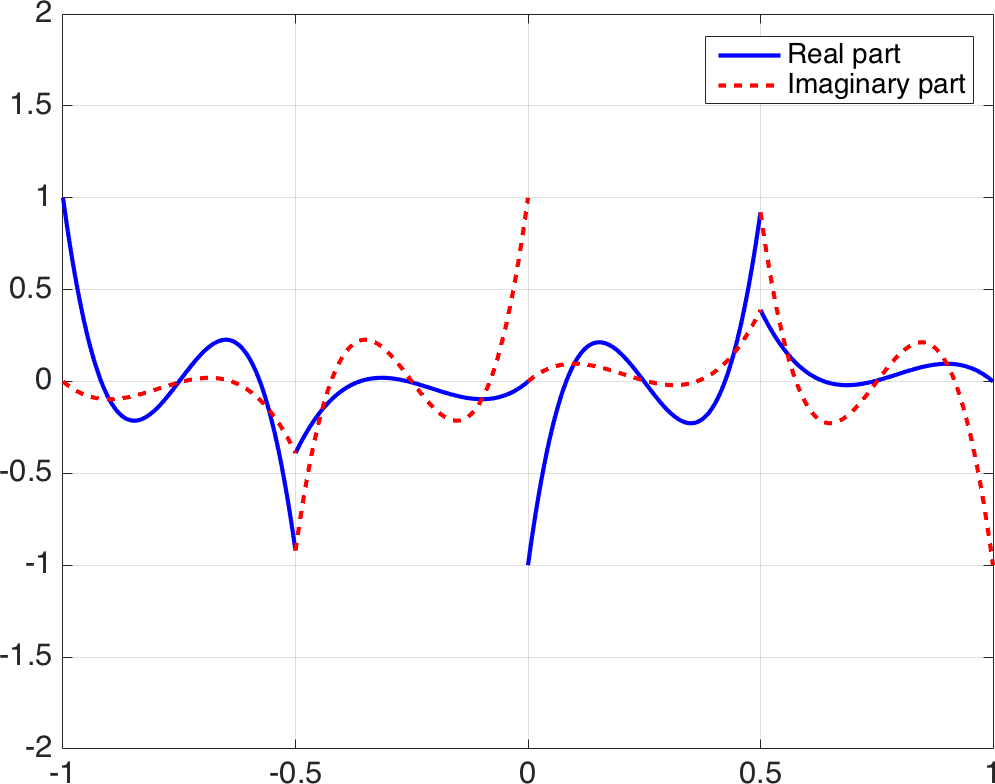}}
\hspace{.5em}
\subfloat[$\tau = 1$]{\includegraphics[width=.32\textwidth]{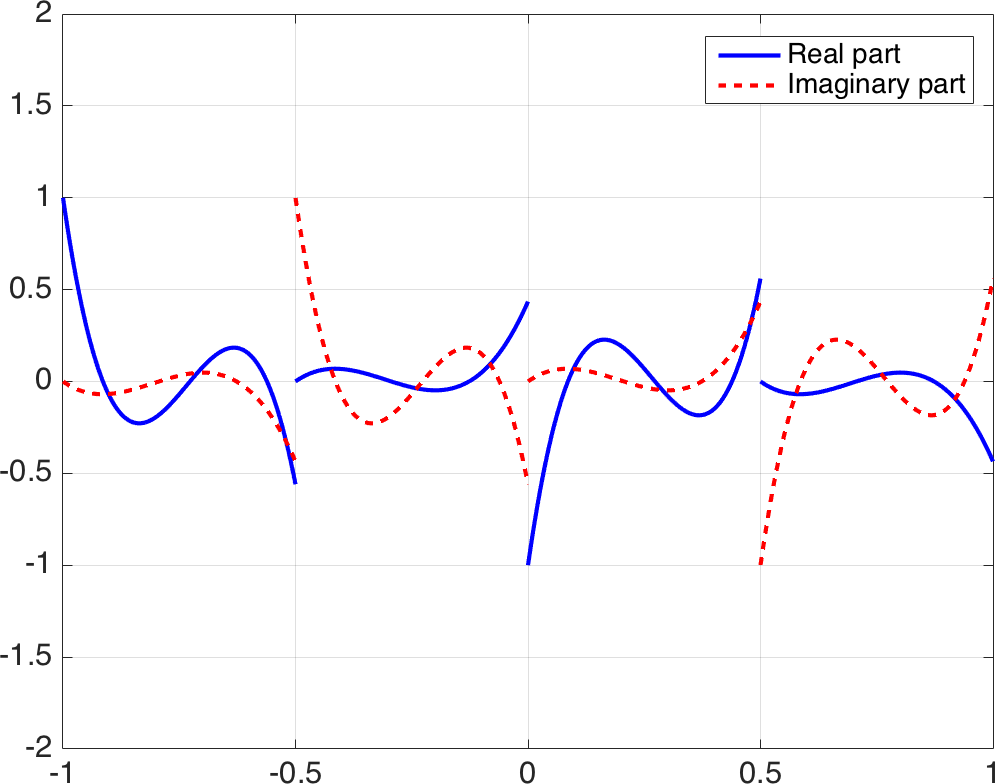}}
\hspace{.5em}
\subfloat[$\tau = 10$]{\includegraphics[width=.32\textwidth]{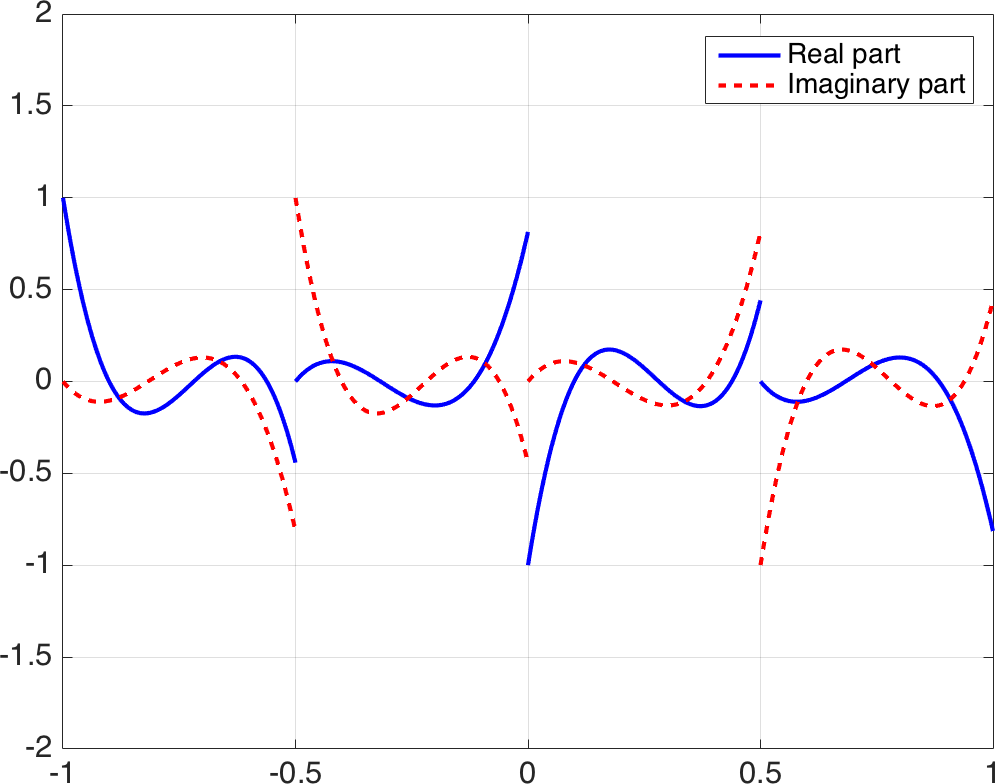}}
\caption{Behavior of a mode for the advection equation corresponding to a divergent eigenvalue (indicated by the red star) as $\tau$ increases.  The order of approximation is $N=3$ on a mesh of 4 elements.} 
\label{fig:trackmodes1D}
\end{figure}

We consider first the scalar advection equation in 1D with $\bm{\beta} = 1$ and periodic boundary conditions.  For $\tau = 0$, the eigenvalues of the DG discretization matrix are purely imaginary.  Paths taken by eigenvalues as $\tau$ increases are determined by sampling \note{spectrum} over a sufficiently fine set of $\tau$ and using a particle tracking method \cite{simpletracker}.  Figure~\ref{fig:track1D} shows the paths taken by these eigenvalues as $\tau$ increases from zero to $\tau = 4$.  As predicted, a subset of \emph{divergent} eigenvalues move further left of the imaginary axis as $\tau$ increases.  The corresponding eigenmodes are shown in Figure~\ref{fig:trackmodes1D}, with inter-element jumps of these eigenfunctions increasing as $\tau$ increases.  

\begin{figure}
\centering
\subfloat{\includegraphics[width=.32\textwidth]{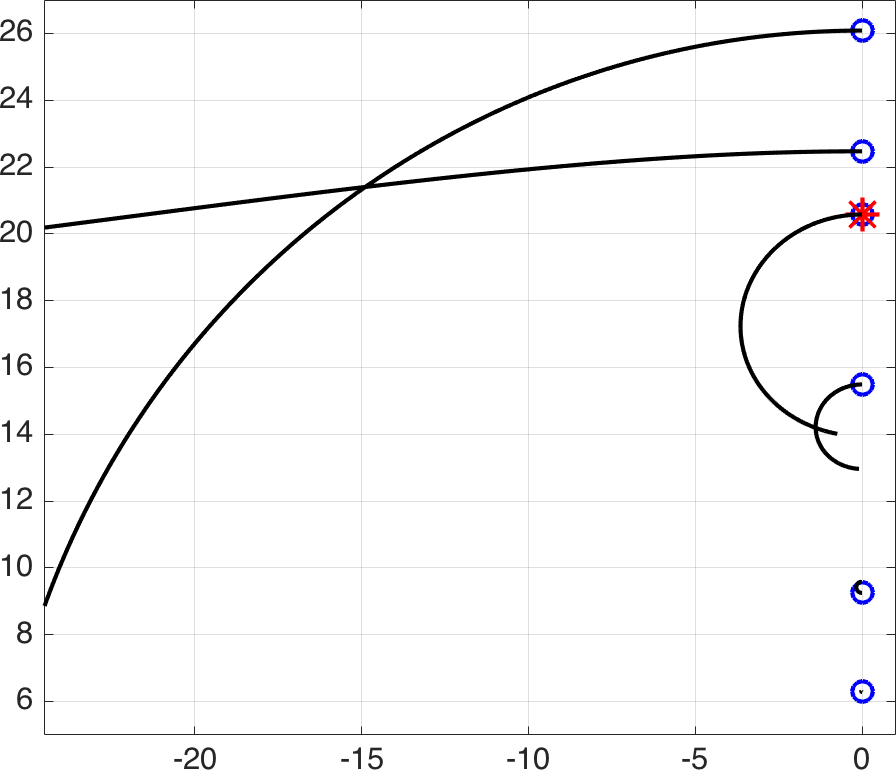}}
\hspace{.5em}
\subfloat{\includegraphics[width=.32\textwidth]{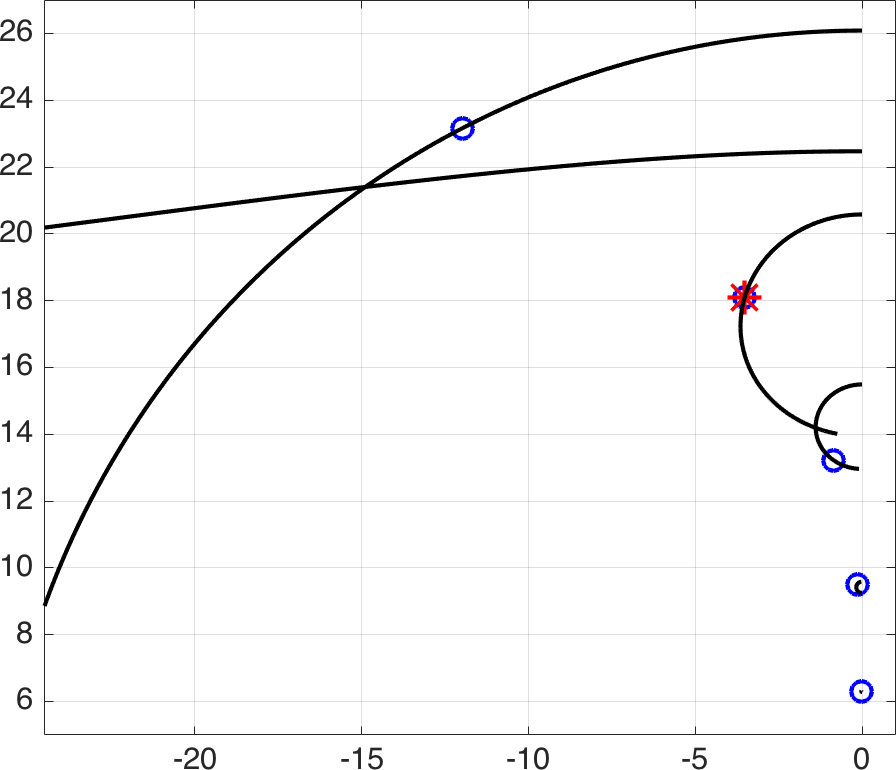}}
\hspace{.5em}
\subfloat{\includegraphics[width=.32\textwidth]{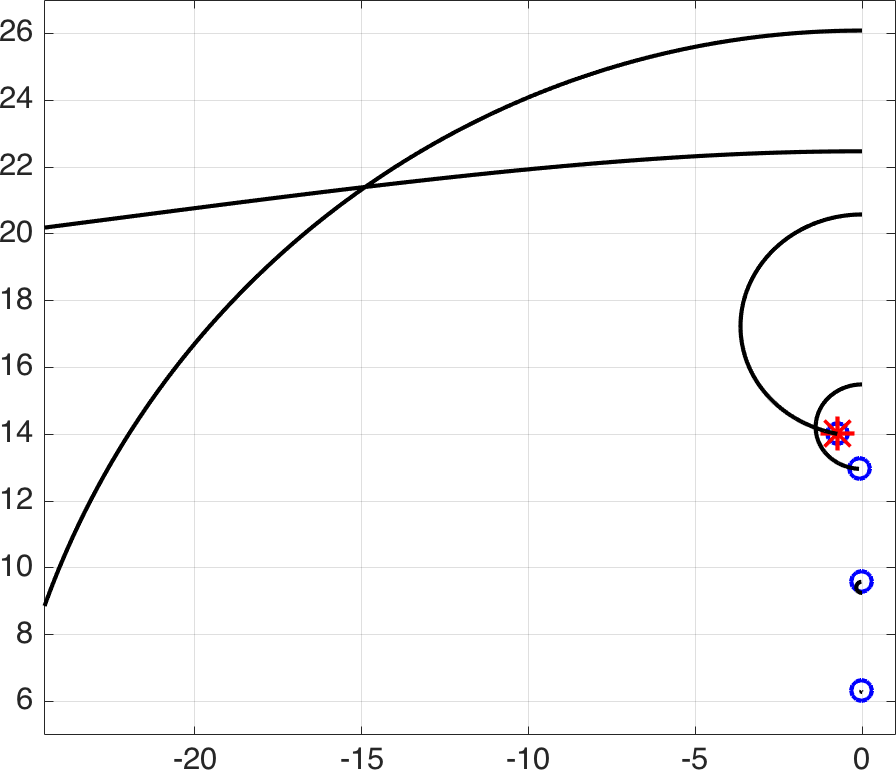}}\\
\subfloat[$\tau = 0$]{\includegraphics[width=.32\textwidth]{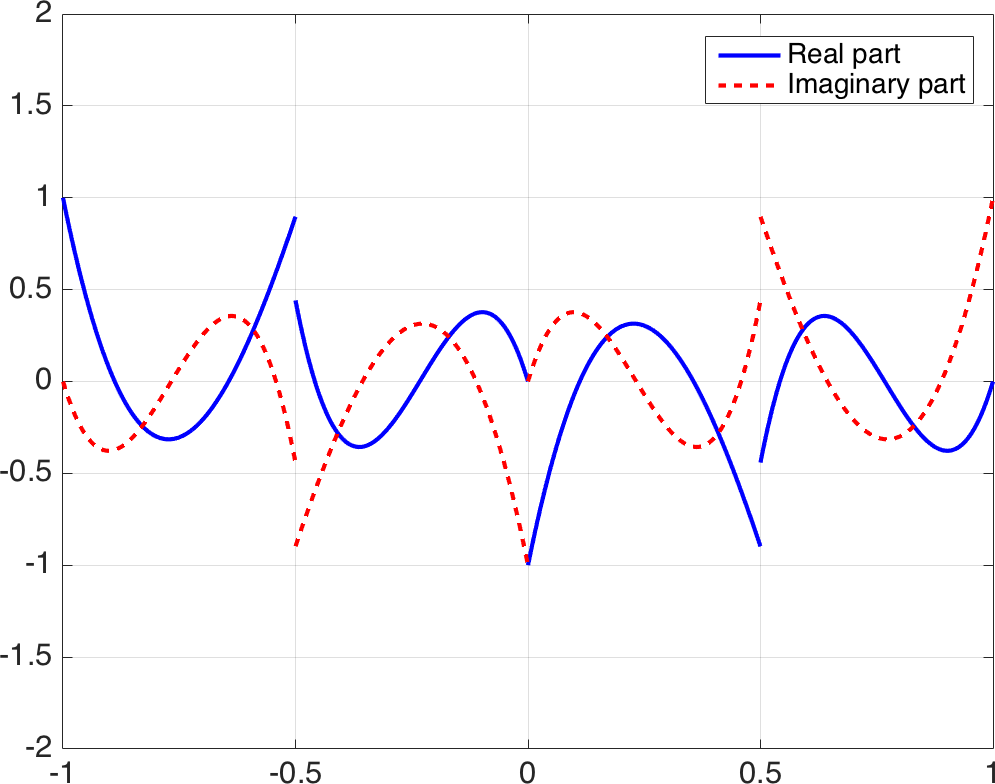}}
\hspace{.5em}
\subfloat[$\tau = 1$]{\includegraphics[width=.32\textwidth]{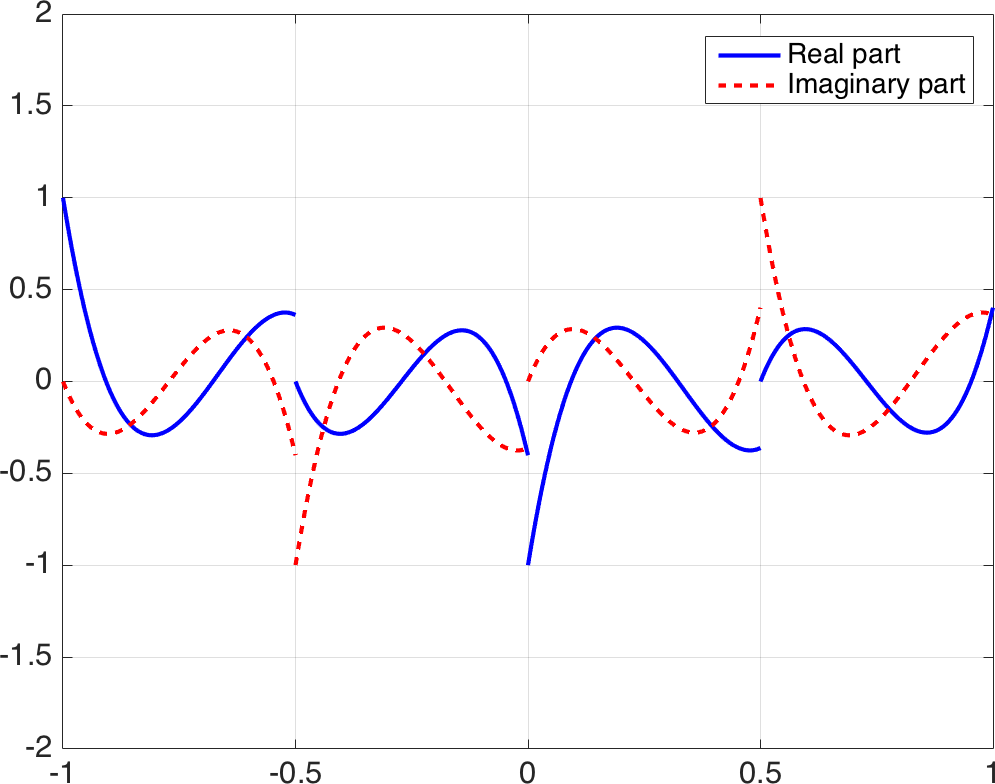}}
\hspace{.5em}
\subfloat[$\tau = 10$]{\includegraphics[width=.32\textwidth]{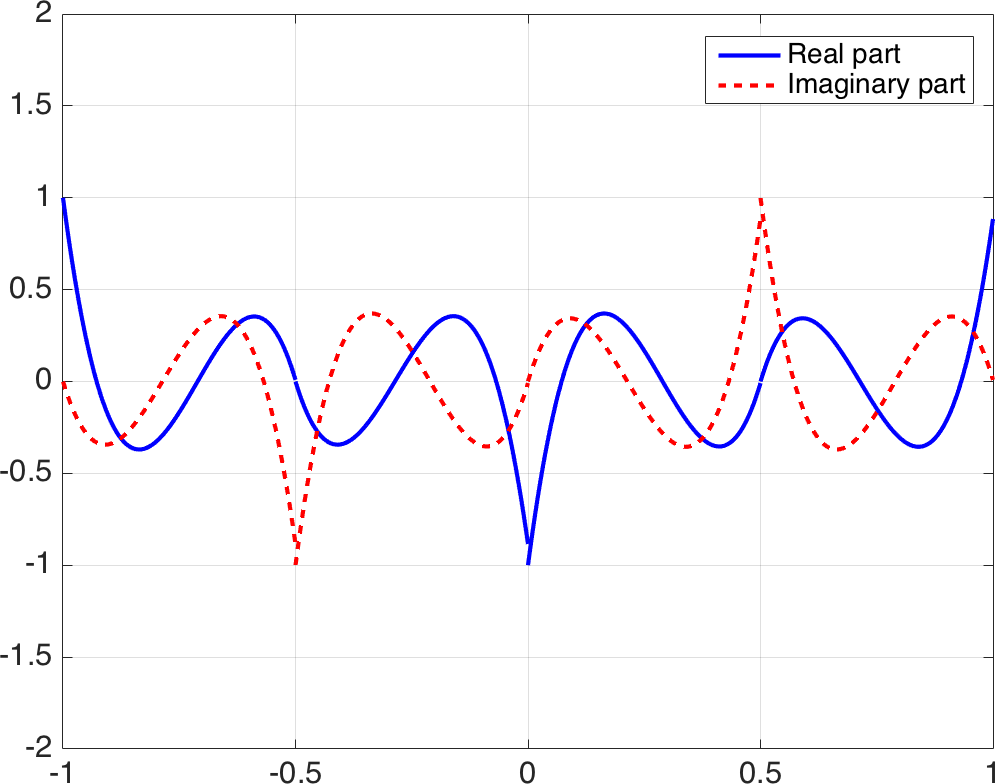}}
\caption{Return behavior of DG advection eigenvalues for sufficiently large $\tau$.  For $\tau = 1$ certain eigenvalues have negative real part, and the corresponding eigenmodes are damped.  For $\tau \gg 1$, certain eigenvalues return to the imaginary axis as spurious modes in conforming discretizations. }
\label{fig:spurious}
\end{figure}

However, for sufficiently large $\tau$, a subset of eigenvalues return to the imaginary axis.  Figure~\ref{fig:spurious} illustrates that, as $\tau$ increases, the magnitude of the inter-element jumps present in these modes decreases, and the mode approaches a $C^0$ continuous function.  These eigenvalues which move right to return to the imaginary axis correspond to a second set of \emph{spurious} modes.  These modes closely resemble spurious modes observed for high order $C^0$ finite element methods \cite{boffi2000problem}, which consist of back-propagating high frequency components with sharp peaks \cite{ainsworth2014dispersive}.  Since the real part of these eigenvalues approaches zero as $\tau\rightarrow \infty$, these modes become undamped and persist as spurious solution components in time-domain simulations \cite{hughes2014finite}.  We note that that a similar phenomena has been observed for the solution of steady state problems in \cite{burman2010interior}, where the authors also point out that large penalty parameters result in decreased accuracy of the DG solution.  

\begin{figure}
\centering
\subfloat[$c_j^0$ (spurious modes for  $\tau = 0$)]{\includegraphics[width=.45\textwidth]{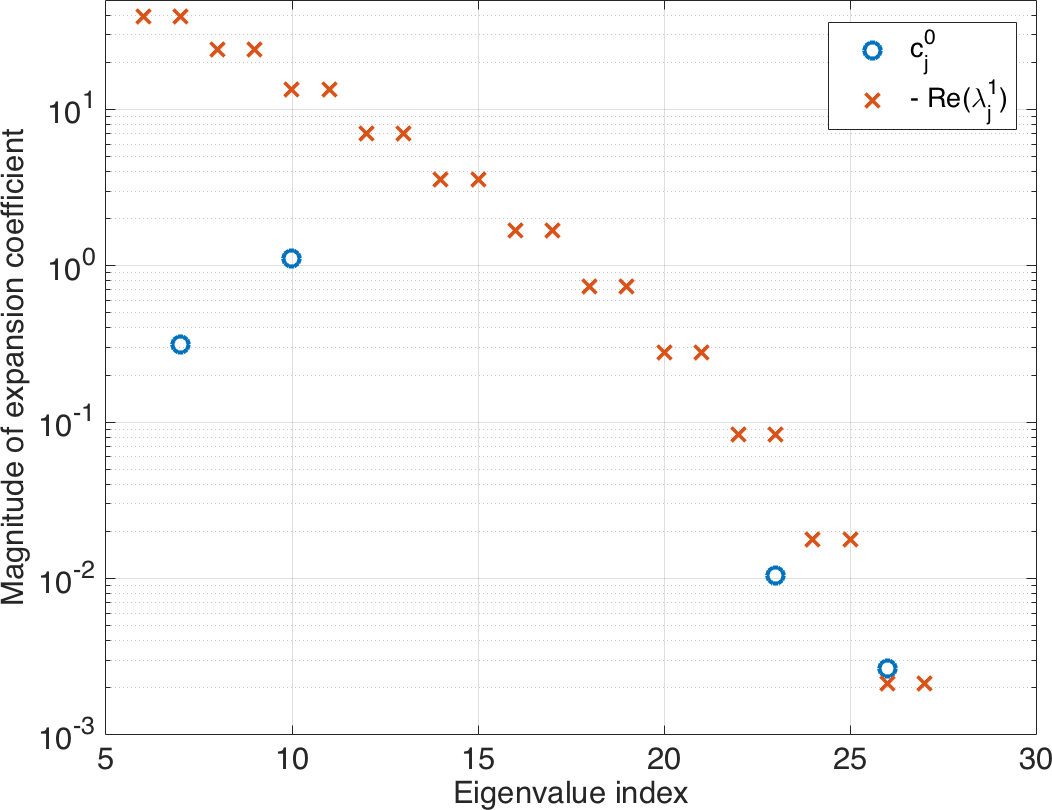}}
\hspace{1em}
\subfloat[$c_j^{100}$ (spurious modes for $\tau = 100$)]{\includegraphics[width=.45\textwidth]{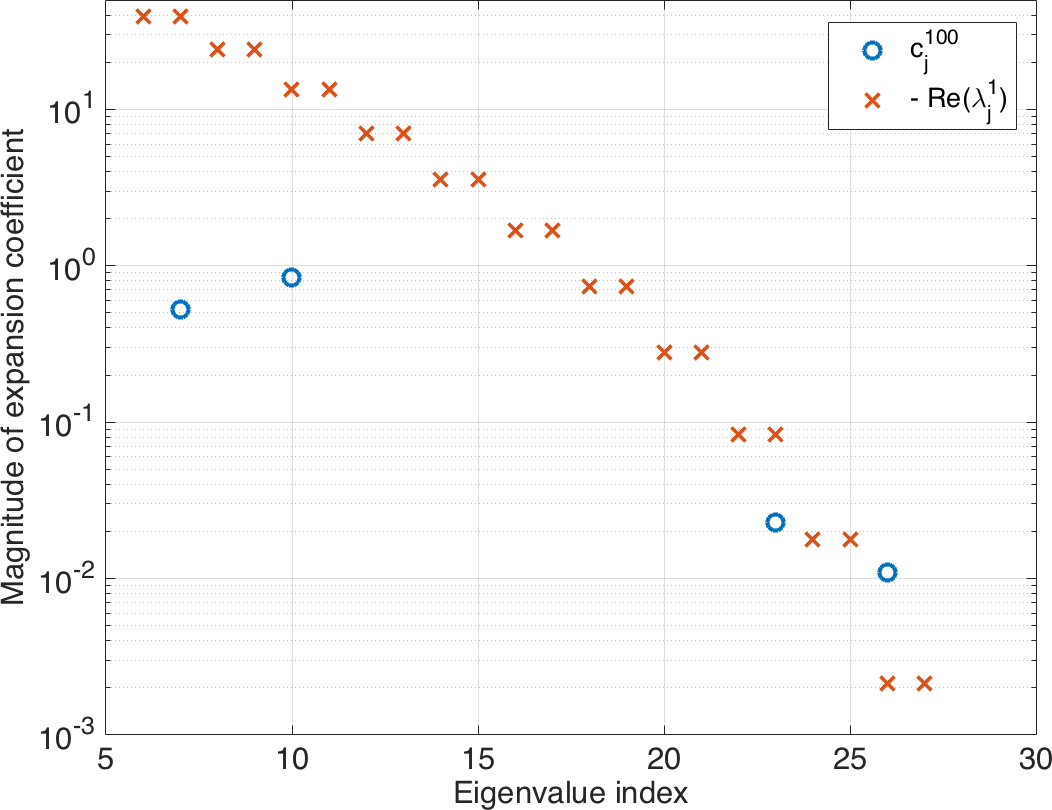}}
\caption{Coefficients $c^0_j$ and $c^{100}_j$ of the spurious mode in Figure~\ref{fig:spurious} when expanded in the modes of $\bm{K}$ for $\tau=1$.  Coefficients with magnitude less than $10^{-13}$ are not shown.   The negative real parts of the corresponding eigenvalues ${\rm Re}(\lambda^1_j)$, which determine the magnitude of damping applied to the $j$th component, are overlaid. }
\label{fig:cmodes}
\end{figure}

Figure~\ref{fig:spurious} also indicates that taking $\tau = 1$ suppresses spurious non-dissipative modes associated with both $\tau = 0$ (central fluxes) and $\tau\rightarrow \infty$ (conforming discretizations).  Let $\lambda^{\tau}_i, \bm{v}^{\tau}_i$ be eigenvalues and eigenvectors of the DG discretization matrix $\bm{K}$ with penalty $\tau$.  The exact solution to the semi-discrete formulation is 
\[
\bm{u}(t) =\bm{V}^{\tau}\LRp{ \sum_j c_j e^{-\lambda_i^{\tau} t} \bm{v}^{\tau}_j}, \qquad \bm{u}(0) = \sum_j c_j\bm{v}^{\tau}_j,
\]
which illustrates how modes $\bm{v}^{\tau}_j$ are damped for $t > 0$ if the corresponding eigenvalue $\lambda^{\tau}_j$ has a negative real part .  We observe that, when spurious modes of the solution for $\tau = 0$ or $\tau \gg 1$ are expanded in the modes for $\tau = 1$
\[
\bm{v}^{\tau}_i = \sum_j c^\tau_j \bm{v}^{1}_j, 
\]
the largest components correspond to eigenvalues with larger negative real parts and heavily damped modes.  As an example, we consider the spurious mode illustrated in Figure~\ref{fig:spurious}.  The absolute values of the coefficients $c^\tau_j$ of this mode are plotted in Figure~\ref{fig:cmodes} for both $\tau = 0$ and $\tau = 100$.  Only the four coefficients larger than $10^{-13}$ are shown.  For comparison, the negative real values of the corresponding eigenvalues $\lambda^1_j$ are also overlaid in Figure~\ref{fig:cmodes}.  The two largest coefficients correspond to eigenvalues with large negative real parts and more strongly damped modes.  The two remaining coefficients correspond to less highly damped eigenvalues with smaller negative real parts, but are an order of magnitude smaller than the two large coefficients.  

The above experiments were repeated for the acoustic wave equation in 1D with reflection boundary conditions on $[-1,1]$, with results similar to those observed for the scalar advection equation.   


\subsection{2D experiments}

In two space dimensions, taking $\tau\rightarrow \infty$ again results in spurious modes which return to the imaginary axis, with the corresponding eigenmodes approaching conforming functions.  We illustrate this through numerical experiments with the acoustic wave equation and advection equation in two space dimensions.  In all numerical experiments, we use $N=3$ and a uniform triangular mesh resulting from bisecting each element of a uniform quadrilateral mesh.  

\begin{figure}
\centering
\subfloat[$\tau = 1$]{\includegraphics[width=.3\textwidth]{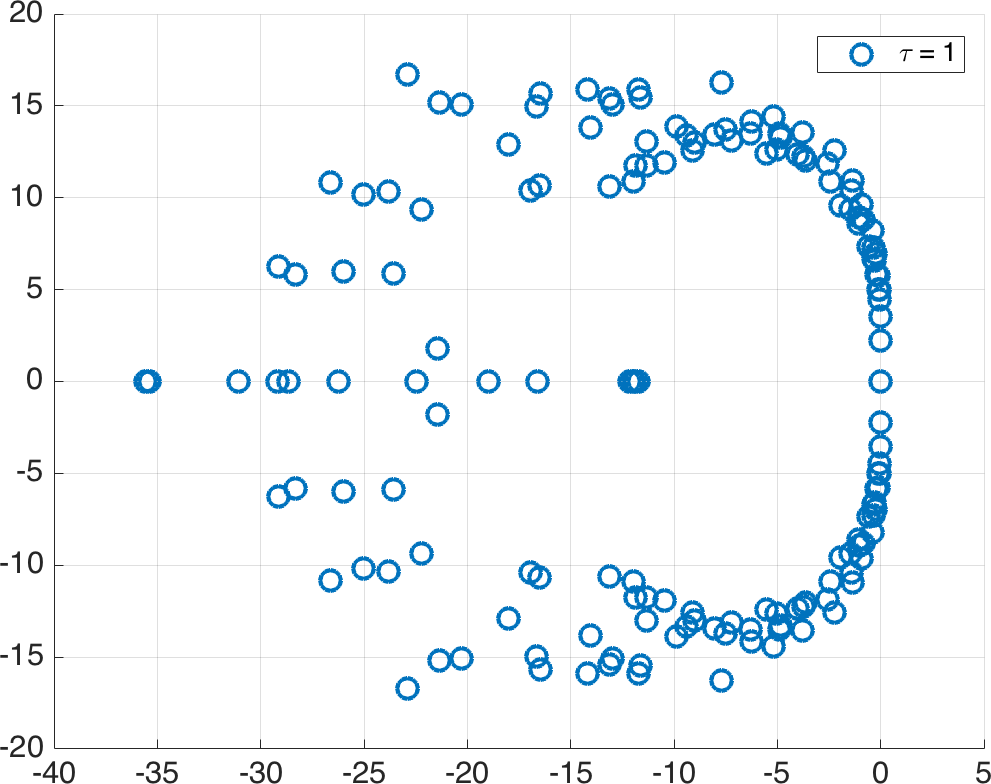}}
\hspace{.5em}
\subfloat[$\tau = 10$]{\includegraphics[width=.3\textwidth]{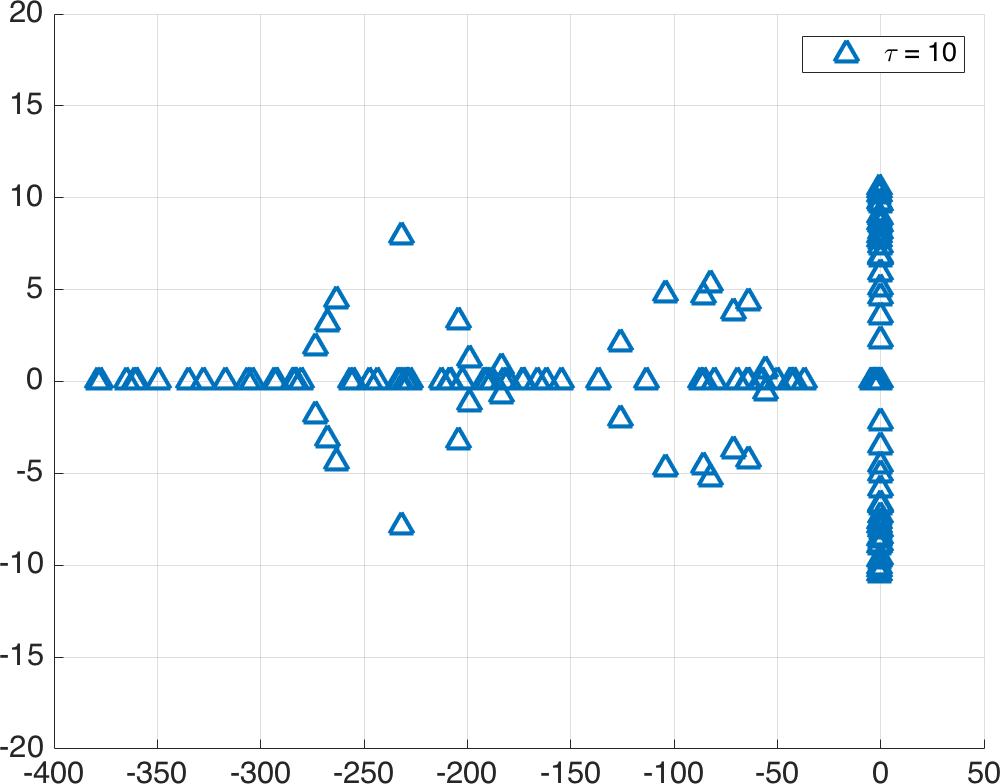}}
\hspace{.5em}
\subfloat[$\tau = 100$]{\includegraphics[width=.3\textwidth]{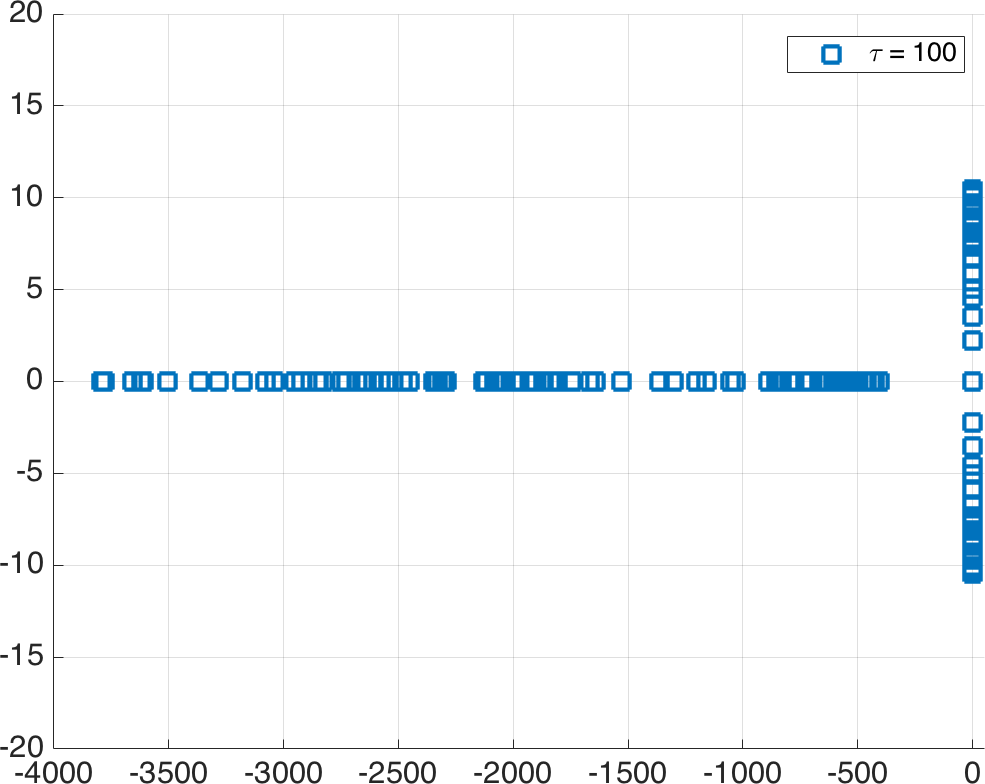}}
\caption{Behavior of eigenvalues for the acoustic wave equation in two dimensions. Note the changing scale of the real axis.}
\label{fig:waveeigs}
\end{figure}

\begin{figure}
\centering
\subfloat[$\tau = 0,1, 50$, zoomed in view]{\includegraphics[height=.3\textheight]{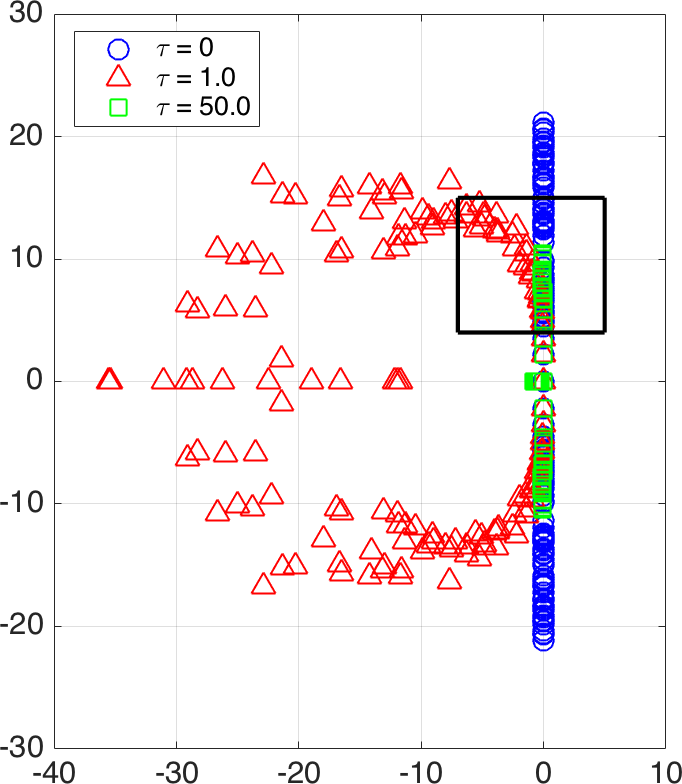}\label{subfig:big}}
\hspace{2em}
\subfloat[$\tau= 0,1,50$, returning eigenvalue paths]{\includegraphics[height=.3\textheight]{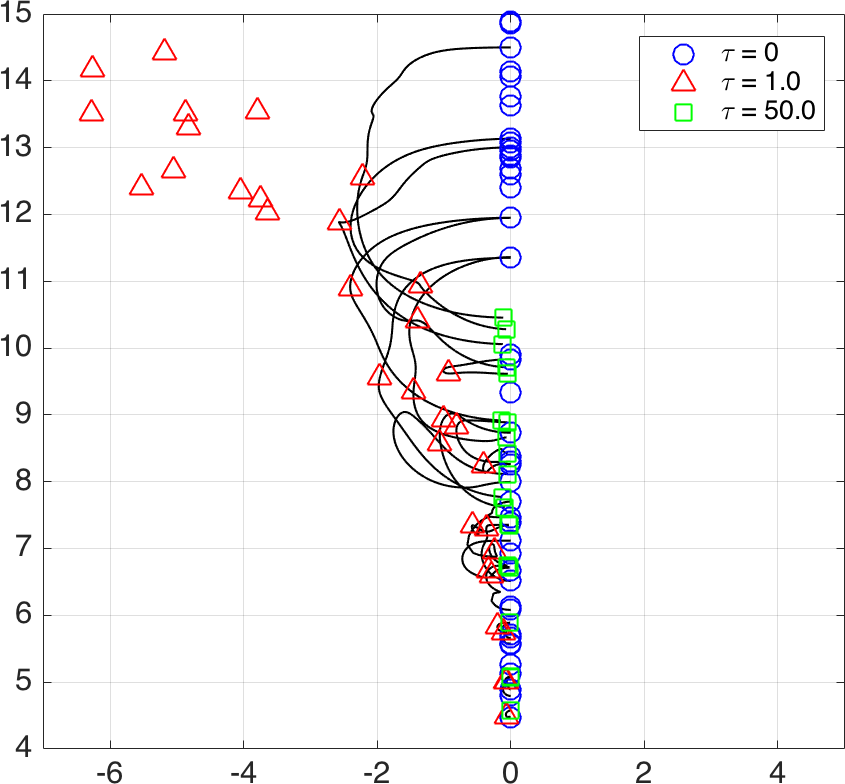}\label{subfig:box}}
\caption{Behavior of eigenvalues for the acoustic wave equation in two dimensions.  Figure~\ref{subfig:box} shows a zoom of the boxed region in Figure~\ref{subfig:big}, with overlaid eigenvalue paths as $\tau$ increases.  Divergent eigenvalues in the far left half plane are not shown.  }
\label{fig:trackEigsWave}
\end{figure}

Figure~\ref{fig:waveeigs} shows the \note{spectrum} of the DG discretization matrix for the acoustic wave equation for $\tau = 1, 10, 50$.  Like the one-dimensional case, a subset of eigenvalues diverge towards the left half plane as $\tau$ increases, though unlike the one-dimensional case, these eigenvalues collapse towards the real line as $\tau$ increases.  Similarly, as predicted, a subset of eigenvalues returns towards the imaginary axis as $\tau \rightarrow \infty$.  These are shown in Figure~\ref{fig:trackEigsWave}, along with traced paths taken by each eigenvalue as $\tau$ increases.  

\begin{figure}
\centering
\subfloat[$\tau = .1$, $\lambda = -0.2379 + 8.7528i$]{\includegraphics[width=.3\textwidth]{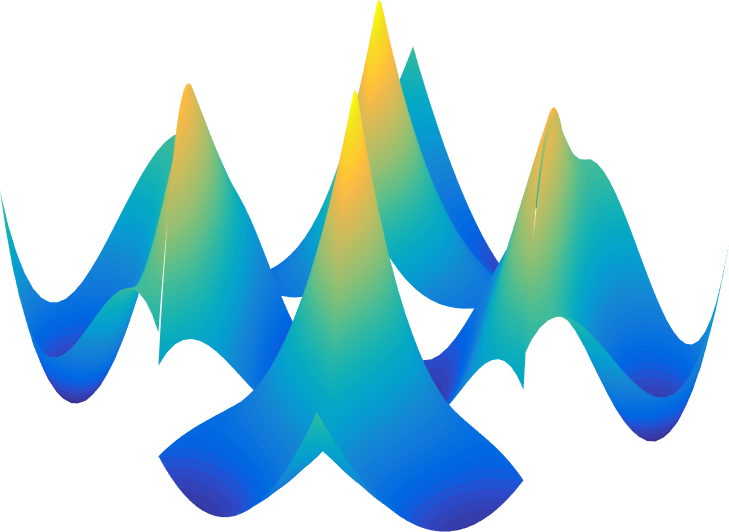}}
\hspace{.5em}
\subfloat[$\tau = 1$, $\lambda = -1.0145 + 8.9270i$]{\includegraphics[width=.3\textwidth]{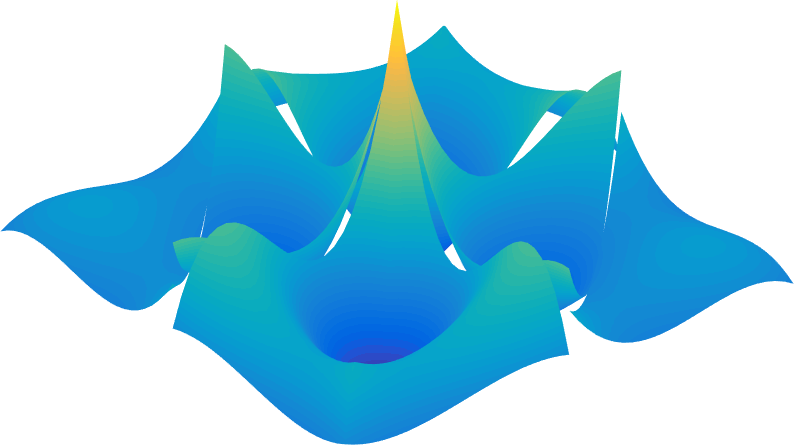}}
\hspace{.5em}
\subfloat[$\tau = 100$, $\lambda = -0.0437 + 7.6167i$]{\includegraphics[width=.3\textwidth]{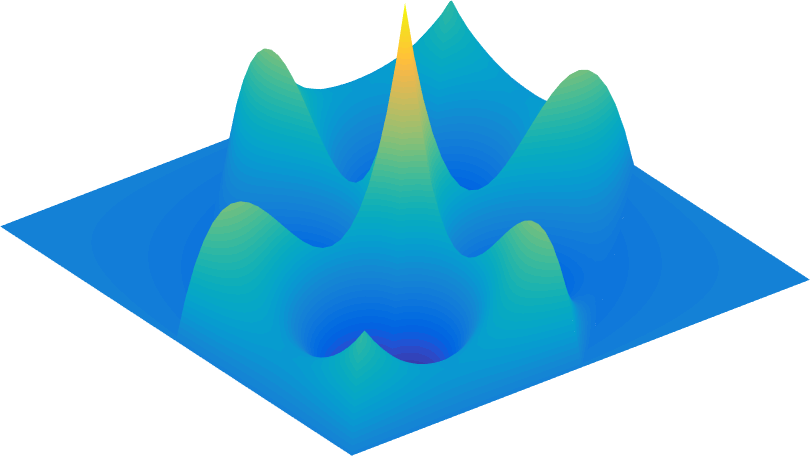}}
\caption{Behavior of the real part of the pressure component of spurious eigenmodes for the acoustic wave equation in two dimensions.  }
\label{fig:trackModesWave}
\end{figure}

The pressure and velocity components of the corresponding acoustic eigenmodes converge to conforming approximations.  As noted in Section~\ref{sec:confexamples}, this implies that pressure components lie in $H^1(\Oh)$ and are continuous across faces, edges, and vertices.  Figure~\ref{fig:trackModesWave} shows the behavior of the pressure component of an eigenmode corresponding to a spurious returning eigenvalue for $\tau = .1, 1, 100$.  As $\tau$ increases and the eigenvalue approaches the imaginary axis, the eigenmode approaches a $C^0$ continuous function with a high frequencies and sharp peak, similar to the spurious eigenmodes observed for the one-dimensional case in Figure~\ref{fig:spurious}.  Likewise, the non-zero velocity components of spurious eigenmodes converge to approximations in $H({\rm div}; \Oh)$, such that only the normal component of velocity is continuous across faces.  


\begin{figure}
\centering
\subfloat{\includegraphics[width=.475\textwidth]{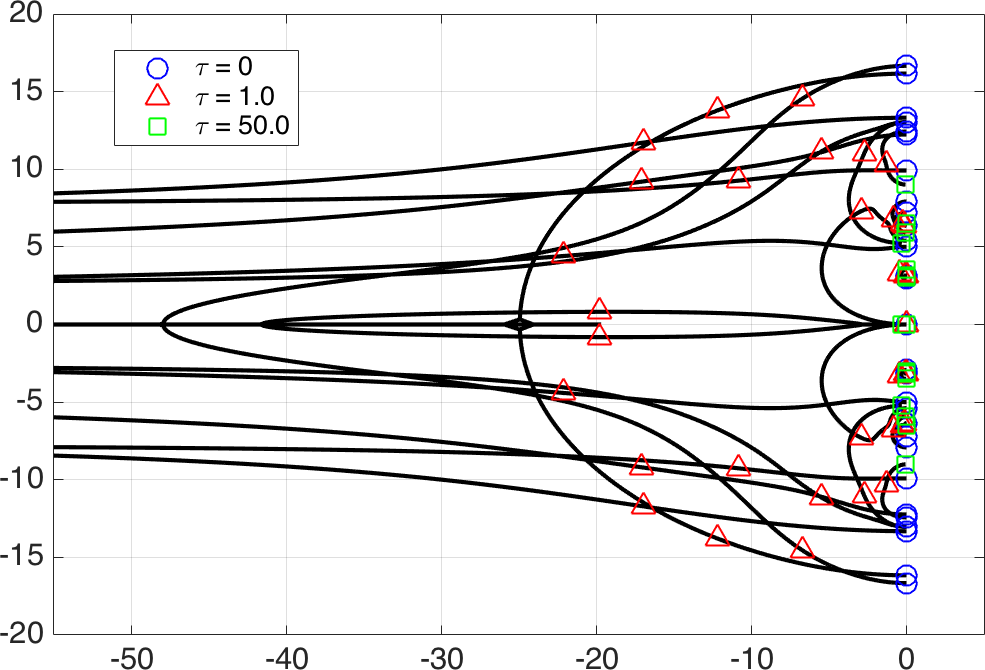}}
\hspace{.5em}
\subfloat{\includegraphics[width=.465\textwidth]{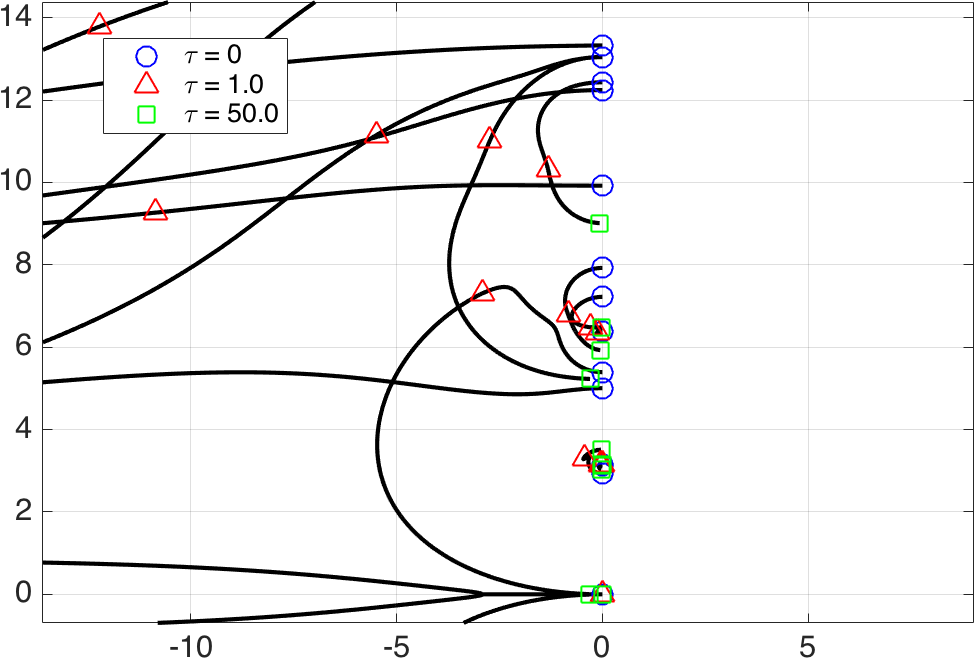}}
\caption{Eigenvalue paths for the DG discretization of advection as $\tau$ increases.  Divergent eigenvalues in the far left half plane are not shown.  }
\label{fig:trackModesAdvec}
\end{figure}

\begin{figure}
\centering
\subfloat[$\tau = 0$, $\lambda = 0 + 7.9246i$]{\includegraphics[width=.3\textwidth]{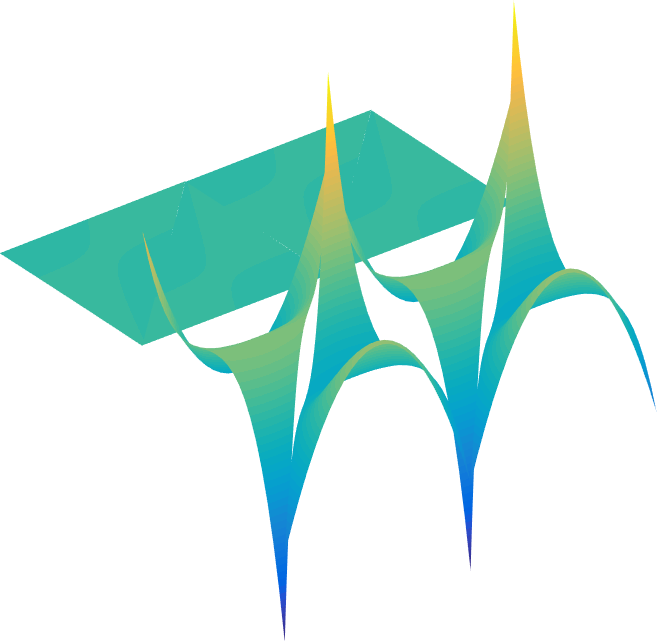}}
\hspace{.5em}
\subfloat[$\tau = 1$, $\lambda =  -0.8278 + 6.7831i$]{\includegraphics[width=.3\textwidth]{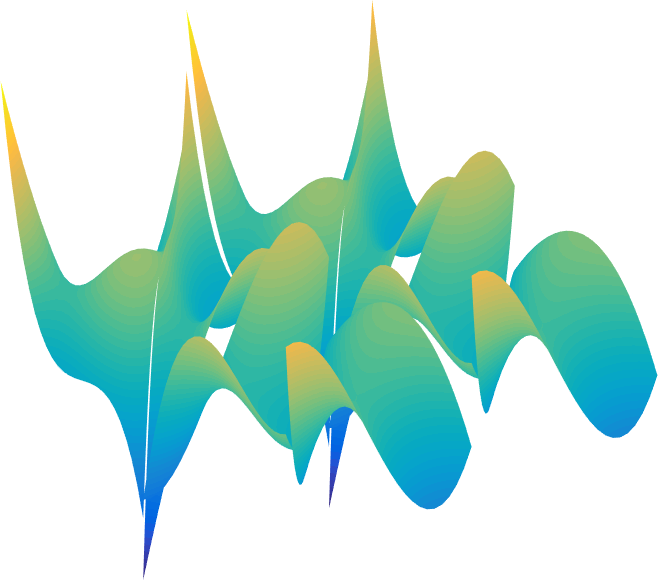}}
\hspace{.5em}
\subfloat[$\tau = 100$, $\lambda = -0.0239 + 5.9282i$]{\includegraphics[width=.3\textwidth]{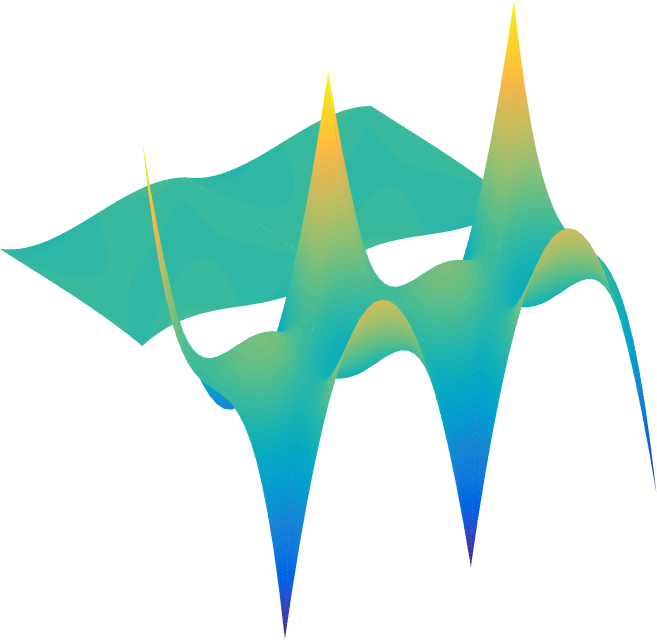}}
\caption{Behavior of the real part of spurious eigenmodes for the advection equation with $\beta = (1,0)$. }
\label{fig:trackModesAdvecU}
\end{figure}

For the advection equation with periodic boundary conditions, the \note{spectrum} of the DG discretization behaves similarly.  In all following experiments, we use advection vector $\beta = (1,0)$.  For $\tau=0$, all eigenvalues lie on the imaginary axis, and as $\tau$ increases, the \note{spectrum} splits into eigenvalues which diverge towards the left half plane and eigenvalues which return to the imaginary axis.  Figure~\ref{fig:trackModesAdvec} shows the \note{spectrum} of the DG matrix for $\tau = 0,1,50$, with eigenvalue paths for increasing values of $\tau$ overlaid.  Unlike spurious conforming modes present for the acoustic wave equation, spurious conforming modes for the advection equation satisfy $\LRb{\bm{\beta}_n}\jump{ u} = 0$, and allow discontinuities along faces which lie tangential to flow directions (Figure~\ref{fig:trackModesAdvecU}).

\section{Conclusions}

For symmetric linear hyperbolic systems of PDEs, penalty fluxes are an alternative to upwind fluxes for DG methods.  Both penalty and upwind fluxes can be interpreted as weakly enforcing continuity conditions implied by the underlying PDE.  As the penalty parameter $\tau$ increases, we show that the eigenvalues of the DG discretization matrix split into two sets, with one set of eigenvalues converging to eigenvalues of a conforming discretization and the other set diverging with real part approaching $-\infty$.  

Numerical experiments also demonstrate that, for $\tau = 1$, there exist eigenvalues with negative real part which return to the imaginary axis as $\tau \rightarrow 0$ or $\tau\rightarrow \infty$.  The corresponding eigenmodes of such eigenvalues can be interpreted as spurious under-resolved modes which persist in time-domain simulations.  These results qualify how the use of upwind and penalty fluxes dampens spurious modes which would otherwise persist for both $\tau = 0$ (central fluxes) and $\tau \rightarrow\infty$ (a conforming discretization).  

\bibliographystyle{unsrt}
\bibliography{dgpenalty}

\end{document}